\documentclass[review,hidelinks,onefignum,onetabnum]{siamart220329}

\usepackage{nicematrix}
\usepackage{tikz,mathrsfs}
\usepackage{amsmath,amsfonts,amssymb}
\usepackage{commath}
\usepackage{mathtools}
\usetikzlibrary{fit}
\usepackage{epsfig}
\usepackage{epstopdf}
\usepackage{float}
\usepackage{hyperref}
\usepackage{comment}
\usepackage{subfigure}
\usepackage{amsfonts}
\usepackage{graphicx}
\usepackage{cite}
\usepackage{textcomp}

\newsiamremark{remark}{Remark}
\newsiamthm{example}{Example}

\nolinenumbers

\title{Partial causal detectability of linear descriptor systems and existence of functional ODE estimators\thanks{Submitted to the editors DATE.
}
}

\author{Juhi Jaiswal\thanks{Department of Chemical Engineering, Indian Institute of Technology Madras, Tamil Nadu - 600036, India (\email{juhijaiswal72@gmail.com}).}
\and Thomas Berger\thanks{Universit\"at Paderborn, Institut f\"ur Mathematik, Warburger Str.~100, 33098~Paderborn, Germany
(\email{thomas.berger@math.upb.de}).}
\and Nutan K. Tomar\thanks{Department of Mathematics, Indian Institute of Technology Patna, Bihar - 801106, India
(\email{nktomar@iitp.ac.in}). This author is the corresponding author.}}

\DeclareMathOperator{\nrank}{\it nor-rank}
\DeclareMathOperator{\row}{row}
\DeclareMathOperator{\esssup}{ess~sup}
\DeclareMathOperator{\blkdiag}{blk-diag}
\DeclareMathOperator{\im}{Im}
\DeclareMathOperator{\rank}{rank}
\DeclareMathOperator{\loc}{loc}

\begin{document}

\maketitle

\begin{abstract}
This paper studies the problem of state estimation for linear time-invariant descriptor systems in their most general form. The estimator is a system of ordinary differential equations (ODEs). We introduce the notion of partial causal detectability and characterize this concept by means of a simple rank criterion involving the system coefficient matrices. Also, several equivalent characterizations for partial causal detectability are established. In addition, we prove that partial causal detectability is equivalent to the existence of functional ODE estimators. A numerical example is given to validate the theoretical results.
\end{abstract}



\section{Introduction}
\noindent We consider linear time-invariant (LTI) descriptor systems of the form
\begin{subequations}\label{dls}
\begin{eqnarray}
E\dot{x}(t) &=& Ax(t)+Bu(t), \label{dlsa} \\
y(t) &=& Cx(t) +Du(t), \label{dlsb} \\
z(t) &=& Kx(t), \label{dlsc}
	\end{eqnarray}
\end{subequations}
where $x : \mathbb{R} \to \mathbb{R}^n,~u: \mathbb{R} \to \mathbb{R}^l,~y: \mathbb{R} \to \mathbb{R}^p$, and $z: \mathbb{R} \to \mathbb{R}^r$ are known as the semistate vector, the input vector, the output vector, and the functional vector, respectively. $E,~A \in \mathbb{R}^{m \times n}$, $B \in \mathbb{R}^{m \times l}$, $C \in \mathbb{R}^{p \times n}$, $D \in \mathbb{R}^{p \times l}$, and $K \in \mathbb{R}^{r \times n}$ with $r \leq n$ are known matrices. The first order matrix polynomial $(\lambda E -A)$, in the indeterminate $\lambda$, is known as matrix pencil. If $m = n$ and $\det(\lambda E -A)$ is a nonzero polynomial in $\lambda$, then system \eqref{dls} is said to be a regular descriptor system. In this article, we consider systems \eqref{dls} in their most general (rectangular) form and assume that the system designer has defined all the coefficient matrices and variables in such a way that the solution set of system \eqref{dls} is non-empty. The tuple $(x,u,y,z):\mathbb{R} \rightarrow \mathbb{R}^{n+l+p+r}$ is said to be a solution of \eqref{dls}, if it belongs to the set
\begin{eqnarray*}
\mathscr{B} &:=& \{(x,u,y,z) \!\in\! \mathscr{L}^1_{\loc}(\mathbb{R}; \mathbb{R}^{n+l+p+r}) \mid Ex \in \mathcal{AC}_{\loc}(\mathbb{R} ; \mathbb{R}^m)     \text{ and } (x,u,y,z)   \text{ satisfies } \\
&&~~~ \eqref{dls} \text{ for almost all } t \in \mathbb{R} \}.
\end{eqnarray*}
Here, $\mathscr{L}^1_{\loc}(\mathbb{R}; \mathbb{R}^{n+l+p+r})$ denotes the set of measurable and locally Lebesgue integrable functions from $\mathbb{R}$ to $\mathbb{R}^{n+l+p+r}$ and $\mathcal{AC}_{\loc} (\mathbb{R} ; \mathbb{R}^m)  $ represents the set of locally absolutely continuous functions from $\mathbb{R}$ to $\mathbb{R}^m$. It is well-known that, corresponding to any given initial condition $Ex(0)$, the system \eqref{dls} may have more than one solution.

In many control applications such as feedback control, fault diagnosis or process monitoring, the information about the full $(K = I_n)$ semistate vector or some part of it is required \cite{trinh2011functional}. However, this information is not available due to physical and/or economical constraints. Hence, in general, the functional vector $z(t) \in \mathbb{R}^r$ contains those variables which cannot be measured and, therefore, we need to estimate them. The existing theory of state estimation for systems of the form \eqref{dls} can be broadly classified in two categories:
\begin{enumerate}
\item[(i)] The estimation generated by a DAE system (described by differential and algebraic equations) of the form
\begin{subequations}\label{obs:dae}
\begin{eqnarray}
E\dot{\hat{x}}(t) &=& A\hat{x}(t)+Bu(t) + L_1 v(t), \label{obs:dae:a} \\
y(t) &=& C\hat{x}(t) +Du(t) + L_2 v(t), \label{obs:dae:b} \\
\hat{z}(t) &=& K\hat{x}(t) + L_3 v(t), \label{obs:dae:c}
\end{eqnarray}
\end{subequations}
where $L_1$, $L_2$,  $L_3$ are matrices of appropriate sizes, and $v(t)$ is an error correction term.
\item[(ii)] The estimation generated by an ODE system (described by ordinary differential equations) of the form
\begin{subequations}\label{obsv}
	\begin{eqnarray}
\dot{w}(t) &=& Nw(t) + H \begin{bmatrix}
u(t) \\ y(t) \end{bmatrix}, \label{obsva} \\
\hat{z}(t) &=& Rw(t) + M \begin{bmatrix}
u(t) \\ y(t) \end{bmatrix}, \label{obsvb}
\end{eqnarray}
\end{subequations}
where $N \in \mathbb{R}^{s \times s}$, $H \in \mathbb{R}^{s \times (l+p)}$, $R\in \mathbb{R}^{r \times s}$, $M \in \mathbb{R}^{r \times (l+p)}$, and $s \in\mathbb{N}\cup\{0\}$.
\end{enumerate}
From an applications point of view, estimation by \eqref{obsv} is always preferred because this system can be initialized arbitrarily and is easily implemented.

In the last few decades, the problem of state estimation for system \eqref{dls} has gained significant attention, due to its wide area of applications in various domains. To the best of our knowledge, the problem of full-state estimation was first considered in $1964$ for state space ($E = I_n$) and in $1983$ for descriptor systems with the seminal works by Luenberger \cite{luenberger1964observing} and El-Tohami et al. \cite{el1983design}, respectively. After this, the theory of full-state estimation for descriptor systems was well developed. Nowadays, there are several equivalent characterizations for the full-state estimation of systems \eqref{dls}, and algorithms for the construction of the estimators exist. A relatively complete literature for the theory of full-state estimation of LTI descriptor systems \eqref{dls} can be found in \cite{berger2017observers,berger2019ode,hou1999observer,jaiswal2021necessary} and the references therein. On the other hand, the problem of functional (or partial-state) estimation has been first addressed in the pioneering work of Dai \cite{dai1989singular} and Minamide et al. \cite{minamide1989design} on regular descriptor systems. In both of these works, the authors estimated $z(t)$ by systems of the form \eqref{obs:dae} under sufficient conditions by fixing $L_2 = I$ and $L_3 = 0$ in system \eqref{obs:dae}. Since then, functional estimators have been used in estimating state space systems with unknown inputs \cite{fernando2007design}, designing observer-based controllers for descriptor systems \cite{ezzine2012controller}, and fault-tolerant controllers for regular descriptor systems \cite{lan2015robust}. In \cite{berger2019disturbance}, Berger studied LTI descriptor systems \eqref{dls} in the context of disturbance decoupled estimation and established a geometric characterization for estimation of the functional vector $z(t)$ via system \eqref{obs:dae}.

Jaiswal et al. \cite{jaiswal2023detectability} introduced the notion of partial detectability for system \eqref{dls} with algebraic as well as geometric characterizations. Further, the authors showed that partial detectability of system \eqref{dls} is necessary for the estimation of the functional vector $z(t)$ via system \eqref{obs:dae}, if $L_2 = I$ and $L_3 = 0$. In this article, we will see that partial detectability is also necessary for the estimation of $z(t)$ via system \eqref{obsv}.

In 2012, Darouach introduced the concept of partial impulse observability as a sufficient condition for the estimation of $z(t)$ \cite{darouach2012functional,darouach2017functional}. Notably, the estimation procedures are correct in \cite{darouach2012functional,darouach2017functional}, but there was a flaw in the algebraic characterization of partial impulse observability. A modified and correct algebraic as well as geometric characterization of partial impulse observability of system \eqref{dls} has been established in \cite{jaiswal2022impulse}. In this article, the authors show that partial impulse observability plays a prominent role in the estimation of the functional vector $z(t)$ by \eqref{obsv}, similar to impulse observability in full-state ($K=I_n$) estimation.

 In $2021$, Jaiswal et al. \cite{jaiswal2021functional,jaiswal2021existence} provided a new set of sufficient conditions for the estimation of $z(t)$ via system \eqref{obsv}, which are weaker than the conditions provided in \cite{darouach2012functional,darouach2017functional,lan2015robust}. In~\cite{jaiswal2024existence}, Jaiswal et al.\ further studied the problem and provided much milder sufficient conditions for the existence of a functional estimator \eqref{obsv}. Although the proposed estimation condition in \cite{jaiswal2024existence} is weaker than all the existing conditions in the literature, it is still not close to being necessary as we will show in Example~\ref{exp2}. In this article, we provide necessary and sufficient conditions for the estimation of~$z(t)$ via systems of the form \eqref{obsv}. 	

The paper is organized as follows. Section \ref{sec:prelim} collects some preliminary results used in the sequel of the article. In Section \ref{sec:Kcausal}, the concept of partial causality  with respect to $K$ for system \eqref{dlsa} is introduced. We provide a rank criterion to test the partial causality of system \eqref{dlsa} with respect to $K$. In addition, this section extends the concept of partial causality to partial causal detectability and establishes several equivalent algebraic and geometric characterizations for the same. In section \ref{sec:ODE:estimator}, necessary and sufficient conditions for the estimation of the functional vector $z(t)$ (via system \eqref{obsv}) are established. A numerical example is given in Section \ref{sec:numerical} to illustrate the step-by-step estimator design procedure. Finally, Section \ref{sec:conc} concludes the article with some future research directions.

We use the following notations throughout the article:  $0$ and $I$ stand for zero and identity matrices of appropriate dimensions, respectively. Sometimes, for more clarity, the identity matrix of size $n \times n$ is denoted by $I_n$.  In a block partitioned matrix, all missing blocks are zero matrices of appropriate dimensions. The set of complex numbers is denoted by  $\mathbb{C}$, $\overline{\mathbb{C}^+} := \{\lambda \in \mathbb{C}~| ~ Re(\lambda)\geq 0\}$ and $\mathbb{C}^- := \{\lambda \in \mathbb{C}~| ~ Re(\lambda)< 0\}$. The symbols $\ker A$, $\row A$, $A^+$, and $A^\top$ denote the null space, the row space, the Moore-Penrose inverse (MP-inverse), and the transpose of a matrix $A \in \mathbb{R}^{m \times n}$, respectively. A matrix pencil $(\lambda E - A)$ is said to have {\it normal rank} $q$ if $\rank (\lambda E - A) = q$ for all but finitely many $\lambda \in \mathbb{C}$ and denoted by $\nrank(\lambda E -A) = q$. In addition, the pencil $(\lambda E - A)$ is
said to be column (row) regular, if it has full column (row) normal rank.  A block diagonal matrix having diagonal elements $A_1,\ldots,A_k$ is represented by $\blkdiag\{A_1,\ldots,A_k\}$. The set $AM := \{Ax~|~x \in M\}$ ($A^{-1}M := \{x \in \mathbb{R}^n~|~Ax \in M\}$) is the image (pre-image) of a subspace $M \subseteq \mathbb{R}^n$ ($M \subseteq \mathbb{R}^m$) under $A \in \mathbb{R}^{m \times n}$. Throughout the article we use the matrices \begin{eqnarray*}
& \mathcal{F}_{n+1,[E,A]} :=
\NiceMatrixOptions
{nullify-dots,code-for-last-col = \color{black},code-for-last-col=\color{black}}
\begin{bNiceMatrix}
E & A &  &   \\
& \Ddots_{(n+1)\text{-times}} & \Ddots &   \\
		&  &  & A  \\
		&  &  & E
	\end{bNiceMatrix} ~\text{ and } ~
	\mathcal{F}_{n+1,[E,A,K]} :=
	\begin{bNiceMatrix}
		E & A &  &  &   \\
		& E & A & & \\
		& &\Ddots_{(n+1)\text{-times}} & \Ddots &    \\
		&  & & & A  \\
		&  &  & & E  \\ & K
	\end{bNiceMatrix}. &
\end{eqnarray*}
For $f\in\mathscr{L}^1_{\loc}(\mathbb{R}; \mathbb{R}^{n})$ we write $f(t) \rightarrow 0$ as $t \rightarrow \infty$, if $\lim\limits_{t \rightarrow \infty} \esssup\limits_{ [t, \infty)} \norm{f(t)}= 0$. 

\section{Preliminaries}\label{sec:prelim}
In this section, we recall some basic concepts from descriptor systems theory and linear algebra. These results will play an important role in the further development of the article.
\begin{lemma}\cite[Quasi-Kronecker Form (QKF)]{berger2013addition}\label{lm:qkf}
For $E,~A \in \mathbb{R}^{m \times n}$ there exist nonsingular matrices $P \in \mathbb{R}^{m \times m}$ and $Q \in \mathbb{R}^{n \times n}$ such that
\begin{equation}\label{eq:qkf}
P(\lambda E - A)Q =  \blkdiag\{  
\lambda E_{\epsilon} - A_{\epsilon}, \lambda I_{n_f} - J_f , \lambda J_{\sigma} - I_{n_\sigma}, \lambda E_{\eta} - A_{\eta}	\},
\end{equation}
where
\begin{enumerate}
\item $E_{\epsilon} ,~ A_{\epsilon} \in \mathbb{R}^{m_{\epsilon} \times n_{\epsilon}}$, $m_{\epsilon} < n_{\epsilon}$, and $\rank (\lambda E_{\epsilon} - A_{\epsilon}) = m_{\epsilon}$, for all $\lambda \in \mathbb{C} \cup \{\infty\}$.
\item $J_f \in \mathbb{R}^{n_f \times n_f}$.
\item $J_{\sigma} \in \mathbb{R}^{n_{\sigma} \times n_{\sigma}}$ is nilpotent.
\item $E_{\eta},A_{\eta} \in \mathbb{R}^{m_{\eta} \times n_{\eta}}$, $m_{\eta} > n_{\eta}$, and $\rank(\lambda E_{\eta} - A_{\eta}) = n_{\eta} $, for all $\lambda \in \mathbb{C} \cup \{\infty\}$.
\end{enumerate}
Here, $\rank (\infty E_{\epsilon} - A_{\epsilon}):=\rank E_{\epsilon}$ and $\rank (\infty E_{\eta} - A_{\eta}):=\rank E_{\eta}$.
\end{lemma}

\begin{remark}\label{rem1}
The blocks in \eqref{eq:qkf} appear only in pairs. For example, if $E_{\epsilon}$ vanishes, then $A_{\epsilon}$ also vanishes. Moreover, $\epsilon-$blocks with $m_\epsilon = 0$ and/or $\eta-$blocks with $n_\eta = 0$ are possible, which results in zero columns (for $m_\epsilon = 0$) and/or zero rows (for $n_\eta = 0$) in the QKF \eqref{eq:qkf}.
\end{remark}

\noindent The following result can be found in any standard textbook of matrix theory.
\begin{proposition}\label{prop:ker}
For matrices $X \in \mathbb{R}^{m \times n}$ and $Y \in \mathbb{R}^{p \times n}$, $\rank \begin{bmatrix}
X \\ Y \end{bmatrix} = \rank X$ if, and only if, $\ker{X} \subseteq \ker{Y}$.
\end{proposition}

The following result is a direct consequence of Proposition \ref{prop:ker}.
\begin{proposition}\label{prop1}
Let $X \in \mathbb{R}^{m \times n}$, $W \in \mathbb{R}^{m \times l}$, $Y \in \mathbb{R}^{p \times n}$, and $Z \in \mathbb{R}^{p \times l}$ be such that
$\rank \begin{bmatrix}
X & W \\ Y & Z \end{bmatrix} = \rank \begin{bmatrix} X & W \end{bmatrix}$, then $\rank \begin{bmatrix}
X \\ Y \end{bmatrix} = \rank{X}$ and $\rank \begin{bmatrix} W \\ Z \end{bmatrix} = \rank{W}$.
\end{proposition}

\begin{proposition}\label{prop2}\cite[Thm.~3.7]{piziak2007matrix}
For matrices $X \in \mathbb{R}^{m \times n}$ and $Y \in \mathbb{R}^{n \times p}$,
\begin{equation*}
\rank(XY) = \rank Y - \dim\left(\ker{X} \cap \im{Y} \right) = \rank X - \dim\left(\ker{Y^\top} \cap \im{X^\top} \right).
\end{equation*}
\end{proposition}

\begin{proposition}\cite{tian2004rank}\label{prop4}
For matrices $X \in \mathbb{R}^{m \times n}$, $W \in \mathbb{R}^{p \times n}$, and $Y \in \mathbb{R}^{p \times l}$,
	$$\rank \begin{bmatrix}
		X & 0 \\ W & Y \end{bmatrix} =  \rank{X} + \rank{Y} + \rank \left((I - YY^+)W(I-X^+X) \right).$$
\end{proposition}
The following result is a simple consequence of Proposition \ref{prop4}
\begin{proposition}\label{prop3}
Let $X \in \mathbb{R}^{m \times n}$, $W \in \mathbb{R}^{m \times l}$, and $Y \in \mathbb{R}^{p \times l}$ be such that $\rank X = m$ and/or $\rank Y = l$, then
$\rank \begin{bmatrix}
X & W \\ 0 & Y
\end{bmatrix} = \rank{X} + \rank{Y}.$
\end{proposition}

\noindent We now recall the following lemma from \cite{jaiswal2021necessary}.
\begin{lemma}\label{lm:Cdecomposition}
Let $E,~A \in \mathbb{R}^{m \times n}$, and $B \in \mathbb{R}^{m \times l}$, then there exist two orthogonal matrices $U_O \in \mathbb{R}^{m \times m}$ and $V_O \in \mathbb{R}^{n \times n}$ such that
\begin{subequations}\label{Cdecomposition}
\begin{eqnarray}
& U_OEV_O =  {\footnotesize \begin{bNiceMatrix}
E_O & E_{k-1} & \boxtimes & \ldots & \boxtimes  \\
& 0 & E_{k-2}  &  \ddots & \vdots  \\
&  & \ddots & \ddots & \boxtimes \\
&  &  & 0 & E_1  \\
& & & & 0  \end{bNiceMatrix}},~U_OAV_O =  {\footnotesize \begin{bNiceMatrix}
A_O & \boxtimes & \hdots & \hdots & \boxtimes   \\
& A_{k-1} & \ddots & & \vdots  \\
&  & \ddots & \ddots & \vdots  \\
&  & & A_2 & \boxtimes  \\
& & &  & A_1  \\
\end{bNiceMatrix}}, & \nonumber \\
&& \label{CdecompositionA} \\
& U_OB = \begin{bNiceMatrix}
	B_O^\top & 0^\top & \ldots & \ldots & 0^\top
	\end{bNiceMatrix}^\top,& \label{CdecompositionB}
		\end{eqnarray}
			\end{subequations}
where $\boxtimes$ represents the matrix elements of no interest and for each $i=1$, $2$, $\hdots$, $k-1$, where $k \leq n$,
\begin{enumerate}
\item[(a)] $A_i$ has full column rank. \label{lm1:a}

\item[(b)] $rank \begin{bmatrix}
\tilde E_i & \tilde B_i
\end{bmatrix} = r_i$, where $r_i$ represents the number of rows in the matrix $\begin{bmatrix}
\tilde E_i & \tilde B_i
\end{bmatrix}$,
$ \tilde E_i = {\footnotesize \begin{bmatrix}
E_O & E_{k-1} & \ldots  & \boxtimes  \\
& \ddots & \ddots & \vdots \\
& & 0 & E_i \end{bmatrix}}$ and $
\tilde B_i=\begin{bmatrix} B_O \\ 0 \end{bmatrix}$.
 \label{lm1:b}

\item[(c)] $\begin{bmatrix}
E_O & B_O \end{bmatrix}$ has full row rank. \label{lm1:c}
	\end{enumerate}
\end{lemma}

The proof of Lemma \ref{lm:Cdecomposition} is given in \cite{jaiswal2021necessary}, and an algorithm to compute $U_O$ and $V_O$ can be found by adapting the similar one in \cite{jaiswal2024existence}. 

Now, we recall the concept of generalized Wong sequences corresponding to a tuple $\{E,A,B,C\}$ from \cite{berger2019disturbance}, various properties of descriptor system \eqref{dls}, and their algebraic and geometric characterizations. It is notable that the original Wong sequences (with $B=0$ and $C=0$) first appeared in a work by Wong \cite{wong1974eigenvalue}, hence their name.

\begin{definition}\label{def:AWS}
For a given system \eqref{dls}, or simply for the tuple $\{E,A,B,C\}$, the generalized Wong sequences
$\left\{\mathcal{V}^i_{[E,A,B,C]} \right\}_{i=0}^\infty$ and  $\left\{\mathcal{W}^i_{[E,A,B,C]} \right\}_{i=0}^\infty$ are sequences of subspaces, defined by
\begin{eqnarray*}%
& \mathcal{V}^{0}_{[E,A,B,C]} := \ker{C}, ~~ \mathcal{V}^{i+1}_{[E,A,B,C]} := A^{-1} (E \mathcal{V}^{i}_{[E,A,B,C]} + \im{B}) \cap \ker{C}, & \\
& \mathcal{W}^{0}_{[E,A,B,C]} := \{0\}, ~~ \mathcal{W}^{i+1}_{[E,A,B,C]} := E^{-1} (A \mathcal{W}^{i}_{[E,A,B,C]} + \im{B}) \cap \ker{C}. &
\end{eqnarray*}
The limits of the generalized Wong sequences are
\begin{eqnarray*} 
& \mathcal{V}^{*}_{[E,A,B,C]} := \bigcap_{i \in \mathbb{N}} \mathcal{V}^{i}_{[E,A,B,C]} ~~ \text{ and }~~ \mathcal{W}^{*}_{[E,A,B,C]} := \bigcup_{i \in \mathbb{N}} \mathcal{W}^{i}_{[E,A,B,C]}.&
\end{eqnarray*}
\end{definition}

\begin{definition}\label{def:completely:controllable}\cite{berger2013controllability}
The descriptor system \eqref{dlsa}, or simply the tuple $\{E,A,B\}$, is completely controllable, if
	\begin{equation*}
		\forall\, x_0,x_f \in \mathbb{R}^n \ \exists\, (x,u,y,z) \in \mathscr{B} \text{ and } t >0 :\ x(0) = x_0 ~\text{ and }~x(t) = x_f.
	\end{equation*}
\end{definition}

\begin{proposition}\cite{berger2013controllability}
The tuple $\{E,A,B\}$ is completely controllable if, and only if, $\mathcal{V}^{*}_{[E,A,B,0]}\cap \mathcal{W}^{*}_{[E,A,B,0]} = \mathbb{R}^n$.
\end{proposition}

\begin{proposition}\label{prop:kcd}\cite{berger2014kalman}
For any $E,~A \in \mathbb{R}^{m \times n}$, $B \in \mathbb{R}^{m \times l}$, and $C\in\mathbb{R}^{p\times n}$, there exist two non-singular matrices $S\in\mathbb{R}^{m\times m}$ and $T\in\mathbb{R}^{n\times n}$ such that
\begin{eqnarray*}
& {\small S\left( \lambda E - A \right)T = \begin{bmatrix}
\lambda E_{11} - A_{11} & \lambda E_{12} - A_{12} & \lambda E_{13} - A_{13} \\
& \lambda E_{22} - A_{22} & \lambda E_{23} - A_{23} \\
&& \lambda E_{33} - A_{33}
\end{bmatrix},
SB = \begin{bmatrix}
			B_1 \\ 0 \\0
		\end{bmatrix},CT = \begin{bmatrix}
		C_1^\top \\ C_2^\top \\ C_3^\top
		\end{bmatrix}^\top,}& 
\end{eqnarray*}
where
\begin{enumerate}
\item[(i)] $E_{11},~A_{11} \in \mathbb{R}^{m_1 \times n_1}$, the triple $\{E_{11},A_{11},B_{1}\}$ is completely controllable, and $m_1 = \rank \begin{bmatrix}
			E_{11} & B_1
\end{bmatrix} \leq n_1 + l$,
		
\item[(ii)] $E_{22},~A_{22} \in \mathbb{R}^{m_2 \times n_2}$ and $E_{22}$ is square ($m_2=n_2$) and invertible,
		
\item[(iii)] $E_{33},~A_{33} \in \mathbb{R}^{m_3 \times n_3}$ with $m_3 \geq n_3$ satisfies $\rank(\lambda E_{33} - A_{33} ) = n_3$ for all $\lambda \in \mathbb{C}$.
\end{enumerate}	
\end{proposition}

We end this section by recalling the concepts of partial impulse observability and partial detectability for system \eqref{dls}. To this end, note that corresponding to inconsistent initial conditions, system \eqref{dls} may possess distributional (impulsive) solutions. Motivated by \cite{berger2017observability}, we denote
\begin{eqnarray*}
\mathscr{B}_{\mathscr{D}} := \{(x,y,z) \in (\mathscr{D}'_{pw \mathscr{C}^{\infty} })^{n+l+p+r} ~|~ (x,y,z) \text{ satisfies \eqref{dls} with $u=0$ on } [0,\infty) \},
\end{eqnarray*}
where $\mathscr{D}'_{pw \mathscr{C}^{\infty}}$ denotes the class of piece-wise smooth distributions and  $\mathscr{B}_{\mathscr{D}}$ is called ITP-behavior in \cite{berger2017observability}. For $f\in \mathscr{D}'_{pw \mathscr{C}^{\infty}}$, the impulsive part at time $t$ is denoted by $f[t]$. For more details, see also \cite{jaiswal2022impulse}.

\begin{definition}\label{def:partial:impulse:observability}\cite{jaiswal2022impulse}
The descriptor system \eqref{dls}, or simply the tuple $\{E,A,C\}$, is said to be partially impulse observable with respect to $K$, if
\[ \forall\, (x,y,z)\in \mathscr{B}_{\mathscr{D}}:\ \big(\forall\, t\ge 0:\ y[t]=0\big)\ \implies\ \big(
\forall\, t\ge 0:\ z[t]=0\big).\]
\end{definition}

In the following lemma, we utilize the fact that
\begin{equation}\label{eq:W*-barEA-EAC}
\mathcal{W}^*_{[\bar E, \bar A, 0, 0]} \cap \ker C = \mathcal{W}^*_{[E,A,0,C]}\text{ and }\bar A^{-1}(\im \bar E) = A^{-1}(\im E) \cap \ker C,
\end{equation}
where the first one follows from Step~4 in the proof of~\cite[Lem.~2.1]{berger2019disturbance} and the second one is clear, and obtain a characterization of partial impulse observability in terms of the generalized Wong sequences.
\begin{proposition}	\cite{jaiswal2022impulse} \label{prop:partial:impulse}
System \eqref{dls} is partially impulse observable with respect to $K$ if, and only if,  $\mathcal{W}^*_{[E,A,0,C]} \cap A^{-1}(\im{E}) \subseteq \ker K. $
\end{proposition}

\begin{definition}\label{def:partial:detectability}\cite{jaiswal2023detectability}
The descriptor system \eqref{dls}, or simply the matrix tuple $\{E,A,C\}$, is said to be partially detectable with respect to $K$, if for all $(x_1,u,y,z_1)$, $(x_2,u,y,z_2)\in \mathscr{B}$ we have
$$z_1(t) - z_2(t) \rightarrow 0 ~\text{ as }~ t \rightarrow \infty.$$
\end{definition}

\begin{proposition}\label{prop:partial:detectability}\cite{jaiswal2023detectability}
The system \eqref{dls} is partially detectable with respect to~$K$ if, and only if, $ \forall\,\lambda \in \overline{\mathbb{C}^+}$,
\begin{equation}\label{eq:partial:detectability}
\rank {\small \begin{bNiceMatrix}
\begin{matrix}
\lambda \bar{E} - \bar{A}
\end{matrix} & & &  \\
\bar{E} &  & &  \\
&\Ddots &  \Ddots^{n \text{- times}} &  \\
&  & \bar{E} & \begin{matrix}
				\lambda \bar{E} - \bar{A}
\end{matrix}   \\
			& & & K
		\end{bNiceMatrix} } = \rank
{\small \begin{bNiceMatrix}
\begin{matrix}
\lambda \bar{E} - \bar{A}
\end{matrix} & & &  \\
\bar{E} &  & &  \\
&\Ddots & \Ddots^{n \text{- times}} &\\
&  & \bar{E} & \begin{matrix}
	\lambda \bar{E} - \bar{A}
	\end{matrix}
\end{bNiceMatrix} }.
\end{equation}
\end{proposition}

\section{Partial causality}\label{sec:Kcausal}
In this section, we first define the concept of partial causality with respect to $K$ for system \eqref{dlsa} and then derive an algebraic criterion in terms of the system coefficient matrices to test the partial causality. The definition of partial causality is a natural extension of {\it causality} of system \eqref{dlsa}, which was introduced by Hou and M\"{u}ller \cite{hou1999causal}. In this section, whenever needed, we take the matrices~$E$ and~$A$ in their QKF \eqref{eq:qkf} to simplify the proofs.
\begin{definition}\label{def:partial:causality}
System \eqref{dlsa}, or simply the triple $\{E,A,B\}$, is said to be  partially causal with respect to $K$, if for every $(x, u, y, z)\in \mathscr{B}$ the system \eqref{dls} has a solution such that $z(t) = Kx(t)$ can be expressed in a form containing no derivatives of $u$.
\end{definition}

In order to analyze this property, we investigate the structure of solutions of \eqref{dls}. The solution theory of descriptor systems is a simple application of the QKF \eqref{eq:qkf} because it has a block diagonal structure and the associated variables can be considered separately. Set
\begin{equation}\label{S2}
x = Q\begin{bmatrix}
x_{\epsilon}^\top & x_f^\top & x_{\sigma}^\top & x_{\eta}^\top
\end{bmatrix}^\top \text{ and }
PB = \begin{bmatrix}
B_{\epsilon}^\top & B_f^\top & B_{\sigma}^\top & B_{\eta}^\top
	\end{bmatrix}^\top ,
\end{equation}
then in terms of the four different blocks in the QKF \eqref{eq:qkf}, \eqref{dlsa} reduces to
\begin{subequations}\label{IV.1}
\begin{eqnarray}
E_{\epsilon}\dot{x}_{\epsilon}(t) &=& A_{\epsilon} x_{\epsilon}(t) + B_{\epsilon}u(t) \label{IV.1a}, \\
\dot{x}_f(t) &=& J_fx_f(t) + B_fu(t) \label{IV.1b},\\
J_{\sigma}\dot{x}_{\sigma}(t) &=& x_{\sigma}(t) + B_{\sigma}u(t) \label{IV.1c}, \\
E_{\eta} \dot{x}_{\eta}(t) &=& A_{\eta}x_{\eta}(t) + B_{\eta} u(t) \label{IV.1d}.
\end{eqnarray}
\end{subequations}
Thus, the following solution analysis of \eqref{dlsa}, via \eqref{IV.1}, is now straightforward. Let $(x,u,y,z) \in \mathscr{B}$ with $x$ partitioned as in \eqref{S2} be given. Then
\begin{enumerate}
\item[S$1$)] in view of assertion $1.$ of Lemma \ref{lm:qkf}, the pencil $(\lambda E_{\epsilon} - A_{\epsilon})$ can (after, possibly, an additional transformation) be written as $\lambda\begin{bmatrix}
	I_{m_\epsilon} & 0
\end{bmatrix} - \begin{bmatrix} A_{\epsilon_1} & A_{\epsilon_2} \end{bmatrix}$.
Therefore, systems of the form \eqref{IV.1a} can also be rewritten as
\begin{equation}\label{S4}
\begin{bmatrix} I_{m_\epsilon} & 0 \end{bmatrix} \begin{bmatrix} \dot x_1 (t) \\ \dot x_2(t) \end{bmatrix} = \begin{bmatrix} A_{\epsilon_1} & A_{\epsilon_2} \end{bmatrix} \begin{bmatrix} x_1(t) \\ x_2(t) \end{bmatrix} + B_{\epsilon}u(t).
\end{equation}
Thus, any solution $x_{\epsilon} = \begin{bmatrix} x_1^\top & x_2^\top \end{bmatrix}^\top$ to \eqref{S4} is given by
\begin{eqnarray*}\label{S1}
\begin{bmatrix}
	x_1(t) \\ x_2(t)
\end{bmatrix} = \begin{bmatrix} \exp(A_{\epsilon_1} t) x_1^0 + \int_0^t \exp(A_{\epsilon_1}(t-\tau)) \left( A_{\epsilon_2} x_2(\tau) + B_{\epsilon}u(\tau) \right) {\rm d}\tau  \\ x_2(t) \end{bmatrix}
\end{eqnarray*}	
for some initial value $x_1^0 \in \mathbb{R}^{m_\epsilon}$. Hence, in general, the system \eqref{IV.1a}  always has a solution, and any solution can be expressed in a form such that $x_{\epsilon}$ contains no derivatives of $u$.
	
\item[S$2$)] Corresponding to any initial condition $x_f^0 \in \mathbb{R}^{n_f}$, the solution of the state space system \eqref{IV.1b} is given by
\begin{equation*}
	x_f(t) = \exp(J_f t) x_f^0 + \int_0^t \exp(J_f(t-\tau)) B_{f} u(\tau) {\rm d}\tau .
\end{equation*}
Therefore, the solution of \eqref{IV.1b} contains no derivatives of $u$.
	
\item[S$3$)] The solution of \eqref{IV.1c} is given by
\begin{equation*}\label{S3}
x_{\sigma}(t) = - \sum_{i = 0}^{h - 1} J_{\sigma}^i B_{\sigma} u^{(i)}(t) ,
\end{equation*}
where $h$ is the nilpotency index of the matrix $J_{\sigma}$, for details see \cite{duan2010analysis}. Hence, the solution of \eqref{IV.1c} contains no derivative of $u$ if, and only if, $u^{(i)}(t) \in \ker(J_{\sigma}^i B_{\sigma})$ for all $0 < i < h$ and for all $t \geq 0$.
	
\item[S$4$)] In view of assertion $4.$ of Lemma \ref{lm:qkf}, the pencil $(\lambda E_{\eta} - A_{\eta})$ can (after, possibly, an additional transformation) be written as $\lambda\begin{bmatrix}
	I_{n_\eta} \\ 0
\end{bmatrix} - \begin{bmatrix} A_{\eta_1} \\ A_{\eta_2} \end{bmatrix}$.
Therefore, systems of the form \eqref{IV.1d} can be rewritten as
\begin{eqnarray*}
\dot{x}_{\eta}(t) &=& A_{\eta_1} x_{\eta}(t) + B_{\eta_1} u(t), \\ 0 &=& A_{\eta_2} x_{\eta}(t) + B_{\eta_2} u(t).
\end{eqnarray*}
Thus, corresponding to any initial condition $x_{\eta}^0 \in \mathbb{R}^{n_\eta}$, the solution is given by
\begin{equation*}
x_{\eta}(t) = \exp(A_{\eta_1} t) x_{\eta}^0 + \int_0^t \exp(A_{\eta_1} (t-\tau))B_{\eta_1} u(\tau) \rm{d} \tau,
\end{equation*}
provided it satisfies $0 = A_{\eta_2} x_{\eta}(t) + B_{\eta_2} u(t)$. Therefore, the solution of \eqref{IV.1d} contains no derivatives of $u$.
\end{enumerate}
In summary, we see that $x$ is forced to contain derivatives of the input $u$ only due to the $\sigma$-block. If the contributions of this block can be excluded from the functional vector~$z$ of system \eqref{dls}, then the system is partially causal. The following result will play an important role in the proof of \Cref{thm:partial:causal:1} below.

\begin{lemma}\label{thm:partial:impulse:free}
Consider system \eqref{dlsa} and \eqref{dlsc}. Then the following statements are equivalent:
\begin{enumerate}
\item $\{E,A,K\}$ satisfies the rank condition \label{lm:Kcausal:c}
\begin{equation}\label{eq:Kcausal:c}
\rank \mathcal{F}_{n+1,[E,A]} =  \rank \mathcal{F}_{n+1,[E,A,K]}.
\end{equation}
		
\item In view of the QKF \eqref{eq:qkf} and $K = \begin{bmatrix}
K_{\epsilon} & K_f & K_{\sigma} & K_{\eta}
\end{bmatrix}$ we have $K_{\sigma}J_{\sigma} = 0$ and $K_{\epsilon} = 0$. \label{lm:Kcausal:a}
		
\end{enumerate}
\end{lemma}

\begin{proof} 
To simplify the rank of $\mathcal{F}_{n+1,[E,A]}$, we apply the following operations:
\begin{enumerate}
\item Write the QKF \eqref{eq:qkf} of $(E,A)$ in each block row.
\item Apply \Cref{prop3} ($(n+1)$-times) from top to bottom to the full row rank matrix $E_{\epsilon}$.
\item Apply \Cref{prop3} ($(n+1)$-times) from top to bottom to the full rank matrix $I_{n_f}$.
\item Apply \Cref{prop3} ($(n+1)$-times) from right to left to the full column rank matrix $E_{\eta}$.
\end{enumerate}
Therefore, we obtain
\begin{eqnarray*}
\rank \mathcal{F}_{n+1,[E,A]} =  (n+1) \left( \rank E_{\epsilon} + \rank I_{n_f} + \rank E_{\eta} \right) + \rank \mathcal{F}_{n+1,[J_{\sigma},I_{n_\sigma}]}.
\end{eqnarray*}

Further, to simplify the rank of $\mathcal{F}_{n+1,[J_{\sigma},I_{n_\sigma}]}$, we apply the following operations:
\begin{enumerate}
\item Multiply $\mathcal{F}_{n+1,[J_{\sigma},I_{n_\sigma}]}$ by $U_{J_{\sigma}} = {\footnotesize \begin{bNiceMatrix}
I_{n_\sigma} &  & & & \\
-J_{\sigma} & \Ddots^{(n+1) \text{- times}} &  & & \\ & \Ddots & & & \\ \Vdots & & & & \\
(-J_{\sigma})^{n} & \Ldots & & -J_{\sigma} & I_{n_\sigma}
\end{bNiceMatrix}}$ from the right.
		
\item Apply \Cref{prop3} ($n$-times) from right to left to the full rank matrix~$I_{n_\sigma}$.
\end{enumerate}
Therefore, in view of the fact that $J_{\sigma}^{n+1} = 0$, we obtain
\begin{eqnarray*}
\rank \mathcal{F}_{n+1,[E,A]} =  (n+1) \left( \rank E_{\epsilon} + \rank I_{n_f} + \rank E_{\eta} \right)  + n \rank I_{n_\sigma}.
\end{eqnarray*}

Similarly, to simplify the rank of $\mathcal{F}_{n+1,[E,A,K]}$, we apply the following operations:
\begin{enumerate}
\item Write the QKF \eqref{eq:qkf} of  $(E,A)$ in the first $(n+1)$-block rows and $K = \begin{bmatrix}
K_{\epsilon} & K_f & K_{\sigma} & K_{\eta}
\end{bmatrix}$ in the $(n+2)^{nd}$-block row.

\item Apply Proposition \ref{prop3} ($(n+1)$-times) from top to bottom to the full rank matrix $I_{n_f}$.

\item Apply Proposition \ref{prop3} ($(n+1)$-times) from right to left to the full column rank matrix $E_{\eta}$.
\end{enumerate}
Therefore, we obtain
\begin{eqnarray*}
\rank \mathcal{F}_{n+1,[E,A,K]} &=& (n+1)(\rank I_{n_f} +  \rank E_{\eta} ) + \\
&& \qquad \rank \mathcal{F}_{n+1,{\tiny \left[\begin{bmatrix}
	E_{\epsilon} \\ & J_{\sigma}
\end{bmatrix},\begin{bmatrix}
A_{\epsilon} \\ & I_{n_\sigma}
\end{bmatrix},\begin{bmatrix}
	K_{\epsilon} & K_{\sigma}
\end{bmatrix}\right]}}  .
\end{eqnarray*}
Additionally, to simplify the rank of $\mathcal{F}_{n+1,{\tiny \left[\begin{bmatrix}
				E_{\epsilon} \\ & J_{\sigma}
			\end{bmatrix},\begin{bmatrix}
				A_{\epsilon} \\ & I_{n_\sigma}
			\end{bmatrix},\begin{bmatrix}
				K_{\epsilon} & K_{\sigma}
			\end{bmatrix}\right]}}$, we apply the following operations:
\begin{enumerate}
\item Multiply $\mathcal{F}_{n+1,{\tiny \left[\begin{bmatrix}
E_{\epsilon} \\ & J_{\sigma}
\end{bmatrix},\begin{bmatrix}
A_{\epsilon} \\ & I_{n_\sigma}
\end{bmatrix},\begin{bmatrix}
K_{\epsilon} & K_{\sigma}
\end{bmatrix}\right]}}$ by
\begin{equation*}
\bar{U}_{J_{\sigma}} = {\footnotesize \begin{bNiceMatrix}
\begin{bmatrix}
I_{n_\epsilon} \\ & I_{n_\sigma}
\end{bmatrix} &  & & & \\
\begin{bmatrix}
	0 \\ & -J_{\sigma}
\end{bmatrix} & \Ddots^{(n+1) \text{- times}} &  & & \\ & \Ddots & & & \\
\Vdots & & & & \\
\begin{bmatrix}
	0 \\ & -J_{\sigma}
\end{bmatrix}^{n} & \Ldots & & \begin{bmatrix}
0 \\ & -J_{\sigma}
\end{bmatrix} & \begin{bmatrix}
	I_{n_\epsilon} \\ & I_{n_\sigma}
\end{bmatrix}
\end{bNiceMatrix}}
\end{equation*} from the right.
		
\item Apply Proposition \ref{prop3} ($n$-times) from right to left to the full rank matrix~$I_{n_\sigma}$.
		
\item Apply Proposition \ref{prop3} to the first block row and full row rank matrix~$E_{\epsilon}$.
\end{enumerate}
Therefore, utilizing the fact that $J_{\sigma}^{n+1} = 0$, we obtain
\begin{eqnarray*}
\rank \mathcal{F}_{n+1,[E,A,K]} &=& \rank E_{\epsilon} + n \rank I_{n_\sigma} + (n+1)(\rank I_{n_f} +  \rank E_{\eta} )  \\
&&~~ + \rank \begin{bNiceMatrix}
\mathcal{F}_{n,[E_{\epsilon},A_{\epsilon}]} \\ \bar{K}_{\epsilon} & K_{\sigma}J_{\sigma}
\end{bNiceMatrix},
\end{eqnarray*}
where $\bar{K}_{\epsilon} = \underset{n \text{-blocks} }{\underbrace{\begin{bmatrix}
K_{\epsilon} & 0 & \ldots & 0
\end{bmatrix}}}$.
Now, by applying Proposition \ref{prop4} to the matrix $\begin{bNiceMatrix}
	\mathcal{F}_{n,[E_{\epsilon},A_{\epsilon}]} \\ \bar{K}_{\epsilon} & K_{\sigma}J_{\sigma}
\end{bNiceMatrix}$ and using the fact that $\rank \mathcal{F}_{n,[E_{\epsilon}, A_{\epsilon}]} = n \rank E_{\epsilon}$, we obtain
\begin{eqnarray*}
 \rank \mathcal{F}_{n+1,[E,A,K]}
&=& (n+1)(\rank E_{\epsilon} + \rank I_{n_f} +  \rank E_{\eta} )  + n \rank I_{n_\sigma} + \rank (K_{\sigma}J_{\sigma}) \\
&& ~ + \rank \left( (I - (K_{\sigma}J_{\sigma})(K_{\sigma}J_{\sigma})^+) \bar{K}_{\epsilon} (I - \mathcal{F}_{n,[E_{\epsilon},A_{\epsilon}]}^+\mathcal{F}_{n,[E_{\epsilon},A_{\epsilon}]}) \right).
\end{eqnarray*}
Thus, rank condition \eqref{eq:Kcausal:c} holds if, and only if,
\begin{eqnarray*}
\rank(K_{\sigma}J_{\sigma}) = 0 &\text{ and }& \rank \left( (I - (K_{\sigma}J_{\sigma})(K_{\sigma}J_{\sigma})^+) \bar{K}_{\epsilon} (I - \mathcal{F}_{n,[E_{\epsilon},A_{\epsilon}]}^+ \mathcal{F}_{n,[E_{\epsilon},A_{\epsilon}]}) \right) = 0 \\
\emph{i.e.},~~ K_{\sigma}J_{\sigma} = 0  &\text{ and }&  \bar{K}_{\epsilon} (I - \mathcal{F}_{n,[E_{\epsilon},A_{\epsilon}]}^+\mathcal{F}_{n,[E_{\epsilon},A_{\epsilon}]}) = 0 \\
\emph{i.e.},~~ K_{\sigma}J_{\sigma} = 0  &\text{ and }& \ker \mathcal{F}_{n,[E_{\epsilon},A_{\epsilon}]} \subseteq \ker \bar{K}_{\epsilon}.
\end{eqnarray*}
We show that $\ker \mathcal{F}_{n,[E_{\epsilon},A_{\epsilon}]} \subseteq \ker \bar{K}_{\epsilon}$ is equivalent to $K_{\epsilon} = 0$. Since $\ker \bar{K}_{\epsilon} = \ker K_{\epsilon} \times \underset{(n-1) \text{-times} }{\underbrace{\mathbb{R}^{n_{\epsilon}} \times \ldots \times \mathbb{R}^{n_{\epsilon}}}}$ it suffices to show  that $\ker \mathcal{F}_{n,[E_{\epsilon},A_{\epsilon}]} \subseteq \ker \bar{K}_{\epsilon}$ implies $K_{\epsilon} = 0$. To this end, let $v_n \in \mathbb{R}^{n_{\epsilon}}$ be arbitrary. Since the Wong sequences terminate after finitely many steps and in each iteration before termination the dimension increases by at least one, we have that $\mathcal{W}^*_{[E_{\epsilon}, A_{\epsilon},0,0]} = \mathcal{W}^{n}_{[E_{\epsilon}, A_{\epsilon},0,0]}$. Furthermore, it is a consequence of~\cite[Lem.~3.11]{BergIlch22} that $\mathcal{W}^*_{[E_{\epsilon}, A_{\epsilon},0,0]} = \mathbb{R}^{n_{\epsilon}}$, thus $v_n\in\mathcal{W}^{n}_{[E_{\epsilon}, A_{\epsilon},0,0]}$. Therefore, there exist $v_i \in \mathcal{W}^{i}_{[E_{\epsilon}, A_{\epsilon},0,0]}$,  $i=1,\ldots,n-1$, such that
\[
 E_{\epsilon}v_{n} + A_{\epsilon}v_{n-1} = 0,~ E_{\epsilon}v_{n-1} + A_{\epsilon}v_{n-2} = 0,~ \ldots,~ E_{\epsilon}v_{2} + A_{\epsilon}v_{1} = 0,~E_{\epsilon}v_1 = 0.
\]
This implies that $\mathcal{F}_{n,[E_{\epsilon}, A_{\epsilon}]}v = 0$ for $v=(v_n^\top,\ldots,v_1^\top)^\top$, hence
\[
    v \in \ker \mathcal{F}_{n,[E_{\epsilon},A_{\epsilon}]} \subseteq \ker \bar{K}_{\epsilon} = \ker K_{\epsilon} \times \underset{(n-1) \text{-times} }{\underbrace{\mathbb{R}^{n_{\epsilon}} \times \ldots \times \mathbb{R}^{n_{\epsilon}}}} \quad\implies\quad v_n \in \ker K_{\epsilon}.
\]
Since $v_n$ was arbitrary, it follows that $\ker K_{\epsilon} = \mathbb{R}^{n_{\epsilon}}$, thus $K_{\epsilon} = 0$. Therefore, we have shown that the rank condition \eqref{eq:Kcausal:c} is equivalent to $K_{\sigma}J_{\sigma} = 0$ and $K_{\epsilon} = 0$.
\end{proof}

The following theorem gives an algebraic characterization of partial causality of system \eqref{dlsa} with respect to $K$, provided $z$ can be  determined uniquely irrespective of~$x$.

\begin{theorem}\label{thm:partial:causal:1}
Consider system \eqref{dlsa}, \eqref{dlsc} and assume that
\begin{equation}\label{eq:assume}
	\nrank \begin{bmatrix} \lambda E - A \\ K \end{bmatrix} = \nrank (\lambda E-A).
\end{equation}
Then the triple $\{E,A,B\}$ is  partially causal with respect to $K$ if, and only if,
\begin{equation}\label{eq:partial:causal:1}
\rank \begin{bNiceMatrix}
\mathcal{F}_{n,[\mathscr{E},\mathscr{A}]} & \mathcal{A} \\  & \mathcal{F}_{n,[E,A]} \end{bNiceMatrix}  = \rank \begin{bNiceMatrix} \mathcal{F}_{n,[\mathscr{E},\mathscr{A}]} & \mathcal{A} \\  & \mathcal{F}_{n,[E,A]} \\ & \mathcal{K} \end{bNiceMatrix},
\end{equation}
where $\mathscr{E}=\begin{bNiceMatrix}
E & 0 \end{bNiceMatrix}$, $\mathscr{A}=\begin{bNiceMatrix}
A & B  \end{bNiceMatrix}$, $\mathcal{A}=\begin{bNiceMatrix} 0 & 0 \\ A & 0 \end{bNiceMatrix}$, and $\mathcal{K} =
\begin{bNiceMatrix}
K & 0 \end{bNiceMatrix}$.
\end{theorem}

\begin{proof}
In view of decomposition \eqref{Cdecomposition} set 	\begin{equation}\label{eq:x:VO}
	x=V_O\begin{bmatrix}
		x_k^\top & x_{k-1}^\top & \ldots & x_1^\top
	\end{bmatrix}^\top \text{ and }~ KV_O = \begin{bmatrix}
	K_O & K_{k-1} & \hdots & K_1
	\end{bmatrix} .
\end{equation}
Also, in view of decomposition \eqref{eq:qkf} of the pencil $(\lambda E_O - A_O)$, set
\begin{equation}\label{eq:x:Q}
x_k = Q\begin{bmatrix}
x_{\epsilon}^\top & x_f^\top & x_{\sigma}^\top & x_{\eta}^\top
\end{bmatrix}^\top \text{ and }~
K_OQ = \begin{bmatrix}
 K_{\epsilon} & K_f & K_{\sigma} & K_{\eta}
\end{bmatrix}.
\end{equation}
Now, we split the proof into the following five steps. \\
\textbf{Step $1$:} In this step, first, we express the assumption \eqref{eq:assume} in terms of the triple $\{E_O,A_O,K_O\}$ and then in terms of the QKF \eqref{eq:qkf} of the pencil $(\lambda E_O - A_O)$. Utilizing decomposition \eqref{Cdecomposition} and \eqref{eq:x:VO} for $K$, as well as assertion $(a)$ of Lemma \ref{lm:Cdecomposition} and Proposition \ref{prop3} for the full column rank matrices $A_i$, $1 \leq i \leq k-1$, we obtain that
\begin{eqnarray}\label{eq:normal:1}
\eqref{eq:assume} \text{ is equivalent to } \nrank \begin{bmatrix} \lambda E_O - A_O \\ K_O \end{bmatrix} = \nrank (\lambda E_O - A_O) .
\end{eqnarray}
Again, by writing the pencil $(\lambda E_O - A_O)$ in the QKF \eqref{eq:qkf}, $K_O$ as in \eqref{eq:x:Q}, and applying Proposition \ref{prop3} for the column regular matrix $\blkdiag\{\lambda I_{n_f} - J_f ,~ \lambda J_{\sigma} - I_{n_\sigma}, ~\lambda E_{\eta} - A_{\eta}\}$, \eqref{eq:normal:1} is equivalent to
\begin{eqnarray}\label{eq:K:epsilon}
 \nrank \begin{bmatrix} \lambda E_{\epsilon} - A_{\epsilon} \\ K_{\epsilon} \end{bmatrix} = \nrank (\lambda E_{\epsilon}-A_{\epsilon}) = m_{\epsilon}, ~ \emph{i.e.},~  K_{\epsilon} = 0
\end{eqnarray}
because the pencil $(\lambda E_{\epsilon} - A_{\epsilon})$ has full row rank for each $\lambda \in \mathbb{C}\cup\{\infty\}$.

\noindent \textbf{Step $2$:} We claim that partial causality of the triple $\{E,A,B\}$ with respect to $K$ is equivalent to partial causality of the triple $\{E_O,A_O,B_O\}$ with respect to $K_O$.

In view of decomposition \eqref{Cdecomposition}, system \eqref{dlsa} and \eqref{dlsc} can be written as
\begin{subequations}\label{dlsdecomposed}
	\begin{eqnarray}
E_O\dot x_k(t) + E_{k-1}\dot x_{k-1}(t)+ \ldots  +\boxtimes\dot  x_1(t) &=& A_Ox_k(t)+ \ldots  + \boxtimes x_1(t) + B_O u(t),\nonumber \\ \label{dlsdecomposedD}\\
E_{k-2}\dot x_{k-2}(t)+\ldots  +\boxtimes\dot  x_1(t) &=& A_{k-1}x_{k-1}(t) + \ldots  +\boxtimes x_1(t), \label{dlsdecomposedC} \\
&\vdots&  \nonumber\\
E_1\dot x_1(t) &=& A_2x_2(t) + \boxtimes x_1(t), \label{dlsdecomposedB}\\
0 &=& A_1x_1(t), \label{dlsdecomposedA} \\
z(t) &=& K_Ox_k(t) + K_{k-1}x_{k-1}(t) + \ldots + K_1x_1(t). \label{dlsdecomposedE}
	\end{eqnarray}
\end{subequations}
Since $A_i$, for $1 \leq i \leq k-1$, has full column rank, solving system \eqref{dlsdecomposed} from \eqref{dlsdecomposedA} to \eqref{dlsdecomposedC}, we obtain
$$x_1 = 0,~x_2 = 0,~\ldots,~ x_{k-1} = 0.$$ Consequently, \eqref{dlsdecomposedD} and \eqref{dlsdecomposedE} reduce to
\begin{subequations}\label{dlsO}
	\begin{eqnarray}
		E_O\dot{x}_k(t) &=& A_Ox_k(t)+B_Ou(t), \label{dlsOa} \\
		z(t) &=& K_Ox_k(t). \label{dlsOb}
	\end{eqnarray}
\end{subequations}
Thus $(x,u,y,z) \in \mathscr{B}$ if, and only if, the tuple $(x_k,u,z)$ satisfies \eqref{dlsO}, where $x = V_O \begin{bmatrix}
	x_k \\ 0
\end{bmatrix}$. This proves the claim.

\noindent \textbf{Step $3$:} In this step, we show that the rank condition \eqref{eq:partial:causal:1} is equivalent to
\begin{equation}\label{eq:partial:causal:2}
\rank \mathcal{F}_{n+1,[E_O,A_O,K_O]} = \rank \mathcal{F}_{n+1,[E_O,A_O]}.
\end{equation}
To simplify the rank of $\begin{bNiceMatrix}
\mathcal{F}_{n,[\mathscr{E},\mathscr{A}]} & \mathcal{A} \\  & \mathcal{F}_{n,[E,A]} \end{bNiceMatrix}$, write the matrices $\mathcal{F}_{n,[\mathscr{E},\mathscr{A}]}$,  $\mathcal{A}$, and $\mathcal{F}_{n,[E,A]}$ in terms of $E,~A,~B$ and substitute decomposition \eqref{Cdecomposition} in the first block row. After that, perform the following operations in the $i^{th}$-row, for $i = 1$ to $i = (k-1)$, repeatedly:
\begin{enumerate}
\item Apply Proposition \ref{prop3} to the full row rank matrix $\begin{bmatrix}
\tilde{E}_i & \tilde{B}_i
\end{bmatrix}$.

\item Substitute decomposition \eqref{Cdecomposition} in the $(i+1)^{st}$-block row.
	
\item Apply Proposition \ref{prop3} to the full column rank matrices $A_j$ in the $i^{th}$-block row, where $1 \leq j \leq i$.
\end{enumerate}
Therefore, we obtain
\allowdisplaybreaks
\begin{eqnarray*}
\rank \begin{bNiceMatrix}
\mathcal{F}_{n,[\mathscr{E},\mathscr{A}]} & \mathcal{A} \\  & \mathcal{F}_{n,[E,A]} \end{bNiceMatrix} &=& \rank {\footnotesize \setcounter{MaxMatrixCols}{25}
\begin{bNiceMatrix}
\tilde{E}_O & \mathscr{A}  &&&&&& & \\ \hline
& \mathscr{E} & \mathscr{A} &&& & &  & \\
&& \Ddots_{(n-k)-\text{blocks}} & \Ddots & & &  & & \\
&&& & \mathscr{A} & & &  &  \\
&&&& \mathscr{E} & A & & & \\ \hline
&&&&& E & A & &  \\
&&&&& & \Ddots_{n-\text{blocks}} & \Ddots &  \\ &&&&& & &  & A  \\  &&&&& & & & E
\end{bNiceMatrix} }   \\
&&~~ +(r_1 + (k-1)\rank A_1)+(r_2+(k-2)\rank A_2)+ \ldots \\
	&&~~ \ldots + (r_{k-1}+\rank A_{k-1}) ,
\end{eqnarray*}
where $\tilde{E}_O = \begin{bmatrix}
	E_O^\top & 0^\top & \ldots & 0^\top
\end{bmatrix}^\top$. Now, substitute $\mathscr{A} = \begin{bmatrix}
	A & B
\end{bmatrix}$, decomposition~\eqref{Cdecomposition} in the $k^{th}$-block row, and perform the following operations in the $i^{th}$-row, for $i = k$ to $i = (n-1)$, repeatedly:
\begin{enumerate}
\item Apply Proposition \ref{prop3} to the full row rank matrix $\begin{bmatrix}
\tilde{E}_O & \tilde{B}_O
\end{bmatrix}$ in the $i^{th}$-block row.
	
\item Substitute decomposition \eqref{Cdecomposition} in the $(i+1)^{th}$-block row.

\item Apply Proposition \ref{prop3} to the full column rank matrices $A_j$ in the $i^{th}$-block row, where $1 \leq j \leq k-1$.
\end{enumerate}
Therefore, we obtain
\begin{align*}
& \rank \begin{bNiceMatrix}
	\mathcal{F}_{n,[\mathscr{E},\mathscr{A}]} & \mathcal{A} \\  & \mathcal{F}_{n,[E,A]} \end{bNiceMatrix} = \rank {\small
	\begin{bNiceMatrix}
	E_O & \tilde{A}_O & & & \\ \hline
	& E & A & &  \\
	& & \Ddots_{n\text{-blocks}} & \Ddots &    \\
	& & & & A \\
	& & & & E
	\end{bNiceMatrix} } + (r_1 + (n-1)\rank A_1) \\
	&\qquad \qquad  +(r_2+(n-2)\rank A_2)+ \ldots  +\big(r_{k-1}+  (n-(k-1)) \rank A_{k-1}\big) \\
	&\qquad \qquad + (n-k)\rank\begin{bmatrix}
		E_O &  B_O
	\end{bmatrix},
\end{align*}
where $\tilde{A}_O = \begin{bmatrix}
	A_O & 0 & \ldots & 0
\end{bmatrix}$. Again, perform the following operations in the $i^{th}$-block row, for $i = n$ to $i = (2n-1)$, repeatedly:
\begin{enumerate}
\item Substitute decomposition \eqref{Cdecomposition} in the $(i+1)^{th}$-block row.
	
\item Apply Proposition \ref{prop3} to the full column rank matrices $A_j$ in the $(i+1)^{th} $-block row, where $1 \leq j \leq k-1$.
\end{enumerate}
Therefore, we obtain
\allowdisplaybreaks
\begin{eqnarray}\label{Frank}
	&& \rank \begin{bNiceMatrix}
		 \mathcal{F}_{n,[\mathscr{E},\mathscr{A}]} & \mathcal{A} \\  & \mathcal{F}_{n,[E,A]} \end{bNiceMatrix} = (r_1 + (n-1)\rank A_1)+(r_2+(n-2)\rank A_2)+ \ldots  \nonumber\\
	&&\qquad\qquad \ldots + \big(r_{k-1} + (n-(k-1))\rank A_{k-1}\big) + (n-k)\rank\begin{bmatrix}
		E_O &  B_O
	\end{bmatrix}  \\
	&& \qquad\qquad+ (n-1)(\rank A_1 +\rank A_2 + \ldots + \rank A_{k-1})	 + \rank \mathcal{F}_{n+1,[E_O,A_O]}. \nonumber
\end{eqnarray}
In a similar manner, we obtain
\begin{eqnarray}\label{F1rank}
	&& \rank \begin{bNiceMatrix} \mathcal{F}_{n,[\mathscr{E},\mathscr{A}]} & \mathcal{A} \\  & \mathcal{F}_{n,[E,A]} \\ & \mathcal{K} \end{bNiceMatrix} = (r_1 + (n-1)\rank A_1)+(r_2+(n-2)\rank A_2)+ \ldots  \nonumber\\
	&&\qquad \qquad \ldots + \big(r_{k-1}+ (n-(k-1))\rank A_{k-1}\big) + (n-k)\rank\begin{bmatrix}
		E_O &  B_O
	\end{bmatrix}  \\
	&&\qquad \qquad + (n-1)(\rank A_1 + \rank A_2 + \ldots + \rank A_{k-1} )	+ \rank \mathcal{F}_{n+1,[E_O,A_O,K_O]}. \nonumber
\end{eqnarray}
Hence, the identities \eqref{Frank} and \eqref{F1rank} reveal that rank condition \eqref{eq:partial:causal:1} is equivalent to
\eqref{eq:partial:causal:2}.

\noindent \textbf{Step $4$:} 	($\Rightarrow$) Assume that \eqref{eq:partial:causal:1} holds. Then, Step $3$ implies that rank condition \eqref{eq:partial:causal:2} holds. Therefore, in view of the QKF \eqref{eq:qkf} for the matrix pencil $(\lambda E_O - A_O)$ and \eqref{eq:x:Q}, Lemma \ref{thm:partial:impulse:free} implies
that $K_{\epsilon} = 0$  and $K_{\sigma} J_{\sigma} = 0$.
Clearly $K_\sigma J_\sigma^i B_\sigma = 0$ for all $i = 1,~2, ~ \ldots,~h-1$, where $h$ is the nilpotency index of $J_{\sigma}$. Therefore, it follows from the solution discussion of \eqref{IV.1} in S$1$)-S$4$) and Definition \ref{def:partial:causality} that the triple $\{E_O,A_O,B_O\}$ is partially causal with respect to $K_O$. Hence, Step~$2$ implies that the triple $\{E,A,B\}$ is partially causal with respect to $K$.

\textbf{Step $5$:} ($\Leftarrow$) Assume that the rank condition \eqref{eq:assume} holds and $\{E,A,B\}$ is partially causal with respect to $K$. Then, Step $2$ implies that $\{E_O,A_O,B_O\}$ is partially causal with respect to $K_O$. By Lemma \ref{lm:Cdecomposition}\,(c), the matrix $\begin{bmatrix} E_O & B_O \end{bmatrix}$ has full row rank. Let $P$ and $Q$ be two nonsingular matrices such that $P(\lambda E_O - A_O)Q$ is in the QKF \eqref{eq:qkf} and $PB_O$ is partitioned as in \eqref{S2}, then
$$P\begin{bmatrix} E_O & B_O	\end{bmatrix}\begin{bmatrix}
	Q \\ & I
\end{bmatrix} = {\footnotesize \begin{bmatrix}
E_{\epsilon} & & & & B_{\epsilon} \\
	& I_{n_f} & & & B_f \\
	& & J_{\sigma} & & B_{\sigma} \\
	& & & E_{\eta} & B_{\eta}
\end{bmatrix}}.$$
By singular value decomposition (SVD) there exist non-singular matrices $U_1$ and $V_1$ such that $E_{\eta} = U_1 \begin{bmatrix} \Sigma_{\eta} \\ 0 \end{bmatrix} V_1^\top$ and $\Sigma_{\eta}$ is invertible. Set $U_2 = \begin{bmatrix}
V_1 \Sigma_{\eta}^{-1} & 0 \\ 0 & I
\end{bmatrix}U_1^\top$. Then $U_2 E_{\eta} = \begin{bmatrix} I_{n_\eta} \\ 0 \end{bmatrix}$ and $U_2 B_{\eta} = \begin{bmatrix} B_{\eta_{11}} \\ B_{\eta_{21}} \end{bmatrix}$. Since $\begin{bmatrix} E_O & B_O \end{bmatrix}$ has full row rank, $B_{\eta_{21}}$ has full row rank as well. Again, it follows from the SVD of $B_{\eta_{21}}$ that there exist non-singular matrices $U_3$ and $V_3$ such that
$ B_{\eta_{21}} = U_3 \begin{bmatrix}
\Sigma_{\eta,2} & 0 \end{bmatrix} V_3^\top$ and $\Sigma_{\eta,2}$ is invertible. Hence, it is clear that there exist invertible matrices $S_1$ and $T_1$ such that
\begin{eqnarray*}
S_1\begin{bmatrix} E_O & B_O \end{bmatrix}T_1 =  {\footnotesize \begin{bmatrix} E_{\epsilon} & & & & B_{\epsilon,1} & 0 \\ & I_{n_f} & & & B_{f,1} & 0 \\ & & J_{\sigma} & & B_{\sigma,1} & 0 \\ & & & I_{n_\eta} & B_{\eta,1} & 0 \\ & & & 0 & 0 & I_{m_\eta - n_\eta}
\end{bmatrix}}.
\end{eqnarray*}
Consequently, invoking the full row rank of $E_\epsilon$, the assumption that the matrix $\begin{bmatrix} E_O & B_O \end{bmatrix}$ has full row rank is equivalent to the fact that $\begin{bmatrix} J_{\sigma} & B_{\sigma,1} \end{bmatrix}$ has full row rank. In view of this decomposition, in the new coordinates the matrix $\begin{bmatrix} \lambda E_O - A_O & B_O \end{bmatrix}$ becomes
\begin{align*}
& S_1 \begin{bmatrix} \lambda E_O - A_O & B_O \end{bmatrix} T_1 \\ & \qquad \qquad  ={\footnotesize\begin{bmatrix} \lambda E_{\epsilon} - A_{\epsilon} & & & & B_{\epsilon,1} & 0 \\ & \lambda I_{n_f} - J_f & & & B_{f,1} & 0 \\ & & \lambda J_{\sigma} - I_{n_\sigma} & & B_{\sigma,1} & 0 \\ & & & \lambda I_{n_\eta} - A_{\eta_1} & B_{\eta,1} & 0 \\ & & & -A_{\eta_2} & 0 & I_{m_\eta-n_\eta}
	\end{bmatrix}}.
\end{align*}
Since the triple $\{E_O,A_O,B_O\}$ is partially causal with respect to $K_O$, it follows from the discussion of the solutions of \eqref{IV.1} in S$1$)-S$4$) (applied to~\eqref{dlsO}) that  $K_\sigma J_\sigma^i B_{\sigma,1} u^{(i)}(t) = 0$ for all $t \geq 0$, $i=1,2,\ldots,h-1$, and for arbitrary $(x,u,y,z) \in \mathscr{B}$. Equivalently,
$K_{\sigma}J_{\sigma}^i B_{\sigma,1} = 0$, for all  $i = 1,~2,~\ldots,~h-1$.
By applying the transposed version of Proposition \ref{prop2} and using the fact that the matrix $\begin{bmatrix} J_{\sigma} & B_{\sigma,1} \end{bmatrix}$ has full row rank, we obtain
\begin{equation*}
	\rank  \left(K_{\sigma}J_{\sigma}^i \begin{bmatrix} J_{\sigma} & B_{\sigma,1} \end{bmatrix}\right) = \rank (K_{\sigma}J_{\sigma}^{i}),~~ \text{ for } i = 1,2,\ldots,h-1 .
\end{equation*}
Thus, for $ 1\leq i \leq h-1$, $K_{\sigma}J_{\sigma}^iB_{\sigma,1} = 0$ implies that $\rank (K_{\sigma}J_{\sigma}^{i+1}) = \rank (K_{\sigma}J_{\sigma}^i)$ and hence
\begin{equation}\label{eq:Sigma}
	\rank (K_{\sigma}J_{\sigma}) = \rank (K_{\sigma}J_{\sigma}^2) = \ldots = \rank (K_{\sigma}J_{\sigma}^{h}) = 0, ~\emph{i.e.,}~ K_{\sigma}J_{\sigma} = 0.
\end{equation}
On the other hand, rank condition \eqref{eq:assume} and Step $1$ imply that $K_{\epsilon} = 0$. Therefore, \eqref{eq:K:epsilon}, \eqref{eq:Sigma}, and Lemma \ref{thm:partial:impulse:free} imply that rank condition \eqref{eq:partial:causal:2} holds. This completes the proof in view of Step $3$.
\end{proof}

\begin{remark}
A careful inspection of the proof of Theorem \ref{thm:partial:causal:1} reveals that the assumption \eqref{eq:assume} is only needed to show that partial causality implies the rank condition \eqref{eq:partial:causal:1}, but not for the converse.
\end{remark}

Now, we extend the definition of {\it partial causality} of \eqref{dlsa} with respect to $K$ to {\it partial causal detectability} of system \eqref{dls} with respect to $K$.

\begin{definition}\label{def:Kcausallydetectable}
System \eqref{dls} is said to be  partially causal detectable with respect to $K$, if the triple  $\{E,A,C\}$ is partially detectable with respect to $K$ and the triple $\{\bar{E},\bar{A},\bar{B}\}$ is partially causal with respect to $K$, where $\bar{E} = \begin{bmatrix}
	E \\ 0 \end{bmatrix}$,
$\bar{A} = \begin{bmatrix}
	A \\ C \end{bmatrix}$, and $\bar{B} = \begin{bmatrix}
	B & 0 \\ D & -I_p \end{bmatrix}$.
\end{definition}

Now, in the following theorem, we derive an algebraic characterization of partial causal detectability with respect to $K$ for system \eqref{dls}.

\begin{theorem}\label{thm:partialcausality}
System \eqref{dls} is partially causal detectable with respect to $K$ if, and only if, the following two rank conditions hold:
\begin{eqnarray}
&\forall\,\lambda \in \overline{\mathbb{C}^+}: ~~\text{ rank condition \eqref{eq:partial:detectability} and} & \label{thm:partialcausalitya} \\
& \rank \begin{bNiceMatrix}
\mathcal{F}_{n,[\mathscr{E},\mathscr{A}]} &  \mathcal{A} \\ & \mathcal{C} \\  & \mathcal{F}_{n,[\bar{E},\bar{A}]} \\ & \mathcal{K} \end{bNiceMatrix} = \rank \begin{bNiceMatrix}
\mathcal{F}_{n,[\mathscr{E},\mathscr{A}]} &  \mathcal{A} \\ & \mathcal{C} \\  & \mathcal{F}_{n,[\bar{E},\bar{A}]}
\end{bNiceMatrix}, & \label{thm:partialcausalityb}
\end{eqnarray}
where  $\mathscr{E} = \begin{bmatrix} E & 0 \end{bmatrix},~\mathscr{A} = \begin{bmatrix}
	A & B \end{bmatrix}$, $ \mathcal{A} = \begin{bmatrix} 0 & 0 \\ A & 0
\end{bmatrix}$, $\mathcal{C} = \begin{bmatrix} C & 0 \end{bmatrix}$, $\bar{E} = \begin{bmatrix}
E \\ 0 \end{bmatrix}$,
$\bar{A} = \begin{bmatrix}
A \\ C \end{bmatrix}$, and $\mathcal{K} = \begin{bmatrix}
K & 0 \end{bmatrix}$.
\end{theorem}

\begin{proof}
($\Rightarrow$): Assume that system \eqref{dls} is partially causal detectable with respect to $K$. Then $\{E,A,C\}$ is partially detectable with respect to $K$ and $\{\bar{E},\bar{A},\bar{B}\}$ is partially causal with respect to $K$. Therefore, it follows from Proposition \ref{prop:partial:detectability} that \eqref{thm:partialcausalitya} holds. Moreover, in view of \Cref{prop1}, condition \eqref{eq:partial:detectability} implies
\begin{equation}\label{eq:unique}
\nrank \begin{bmatrix}
	\lambda \bar{E} - \bar{A} \\ K
		\end{bmatrix} = \nrank (\lambda \bar{E} - \bar{A}).
\end{equation}
Hence, it follows from \eqref{eq:unique}, partial causality of $\{\bar{E},\bar{A},\bar{B}\}$ and Theorem \ref{thm:partial:causal:1} that
\begin{equation}\label{eq:Fbar}
\rank \begin{bNiceMatrix} \mathcal{F}_{n,[\bar{\mathscr{E}},\bar{\mathscr{A}}]} & \bar{ \mathcal{A}} \\  & \mathcal{F}_{n,[\bar{E},\bar{A}]} \end{bNiceMatrix}  = \rank \begin{bNiceMatrix}  \mathcal{F}_{n,[\bar{\mathscr{E}},\bar{\mathscr{A}}]} & \bar{ \mathcal{A}} \\  & \mathcal{F}_{n,[\bar{E},\bar{A}]} \\ & \mathcal{K} \end{bNiceMatrix}.
\end{equation}	
Now, by writing the matrix $\mathcal{F}_{n,[\bar{\mathscr{E}},\bar{\mathscr{A}}]}$ in terms of the system coefficient matrices $E,~A,~B,~C,~D$, and $I_p$, it is easy to see that the identity matrix $I_p$ appears in $(n-1)$ columns corresponding to $\bar{B}$. By permuting these identity matrices to the left upper corner in diagonal positions and applying Proposition \ref{prop2}, we obtain
\begin{eqnarray}
\rank \begin{bNiceMatrix}
\mathcal{F}_{n,[\bar{\mathscr{E}},\bar{\mathscr{A}}]} & \bar{ \mathcal{A}} \\  & \mathcal{F}_{n,[\bar{E},\bar{A}]}
\end{bNiceMatrix}  &=&  \rank \begin{bNiceMatrix}
\mathcal{F}_{n,[\mathscr{E},\mathscr{A}]} &  \mathcal{A} \\ & \mathcal{C} \\  & \mathcal{F}_{n,[\bar{E},\bar{A}]} \end{bNiceMatrix} + (n-1)p , \label{eq:barFn} \\
\rank \begin{bNiceMatrix}
\mathcal{F}_{n,[\bar{\mathscr{E}},\bar{\mathscr{A}}]} & \bar{ \mathcal{A}} \\  & \mathcal{F}_{n,[\bar{E},\bar{A}]} \\ & \mathcal{K} \end{bNiceMatrix} &=& \rank \begin{bNiceMatrix}
\mathcal{F}_{n,[\mathscr{E},\mathscr{A}]} &  \mathcal{A} \\ & \mathcal{C} \\  & \mathcal{F}_{n,[\bar{E},\bar{A}]} \\ & \mathcal{K} \end{bNiceMatrix} + (n-1)p. \label{eq:barFn1}
\end{eqnarray}
Then, Eqs. \eqref{eq:Fbar}, \eqref{eq:barFn}, and \eqref{eq:barFn1} imply rank condition \eqref{thm:partialcausalityb}.
	
($\Leftarrow$): Clearly, condition \eqref{thm:partialcausalitya} implies partial detectability of $\{E,A,C\}$ with respect to $K$ and~\eqref{eq:unique}. In addition, the assumption \eqref{thm:partialcausalityb}, rank identity \eqref{eq:barFn} and \eqref{eq:barFn1} imply that \eqref{eq:Fbar} holds. Therefore, it follows from Theorem \ref{thm:partial:causal:1} that $\{\bar{E},\bar{A},\bar{B}\}$ is partially causal with respect to $K$. This completes the proof.
\end{proof}

By Theorem \ref{thm:partialcausality}, partial causal detectability is characterized by the rank condition \eqref{eq:partial:detectability} for partial detectability together with the rank condition \eqref{thm:partialcausalityb}. The latter is amenable to a variety of further characterizations, which can be found in the following theorem.

\begin{theorem}\label{thm:Kcausallydetectable}
For system \eqref{dls}, the following statements are equivalent:
	\begin{enumerate}
\item[(i)] rank condition \eqref{thm:partialcausalityb} holds. \label{b}
		
\item[(ii)] $\mathcal{A}^{-1}\left(\im{ \mathcal{F}_{n,[\mathscr{E},\mathscr{A}]}}\right) \cap \ker{\mathcal{C}} \cap \ker \mathcal{F}_{n,[\bar{E},\bar{A}]}  \subseteq \ker{\mathcal{K}}$. \label{e}
		
\item[(iii)] $\mathcal{A}_1^{-1}(\im{ \mathcal{F}_{n,[\mathscr{E},\mathscr{A}]}}) \cap \mathcal{W}^{*}_{[E,A,0,C]} \subseteq \ker{K}$, where
$\mathcal{A}_1 = \begin{bmatrix}
	0 \\ A
\end{bmatrix}$. \label{a}

\item[(iv)] $A^{-1}\left(E \left( \mathcal{V}_{[E,A,B,0]}^{n-1} \right) \right) \cap \mathcal{W}_{[E,A,0,C]}^{*} \subseteq \ker K. $ \label{c}
		
\item[(v)] The completely controllable part of system \eqref{dls} is partially impulse observable with respect to the corresponding part of $K$ according to Kalman controllability decomposition from Proposition \ref{prop:kcd}. \label{d}
\end{enumerate}
\end{theorem}

	\begin{proof}
$(i) \Leftrightarrow (ii)$:
Let $\mathcal{Z}$ be any matrix such that $\ker{\mathcal{Z}} = \im{\mathcal{F}_{n,[\mathscr{E},\mathscr{A}]}}$. Then, in view of Proposition \ref{prop2}, we obtain
\begin{eqnarray*}
\rank \begin{bmatrix}
\mathcal{Z}\mathcal{A} \\ \mathcal{C} \\ \mathcal{F}_{n,[\bar{E},\bar{A}]}
\end{bmatrix}
&=& \rank \begin{bNiceMatrix}
\mathcal{F}_{n,[\mathscr{E},\mathscr{A}]} & \mathcal{A} \\ & \mathcal{C} \\ & \mathcal{F}_{n,[\bar{E},\bar{A}]} \end{bNiceMatrix} - \rank
\mathcal{F}_{n,[\mathscr{E},\mathscr{A}]} , \\
\rank \begin{bmatrix} \mathcal{Z}\mathcal{A} \\ \mathcal{C} \\  \mathcal{F}_{n,[\bar{E},\bar{A}]} \\ \mathcal{K} \end{bmatrix} &=& \rank \begin{bNiceMatrix}
\mathcal{F}_{n,[\mathscr{E},\mathscr{A}]} & \mathcal{A} \\ & \mathcal{C} \\
& \mathcal{F}_{n,[\bar{E},\bar{A}]} \\
& \mathcal{K}  \end{bNiceMatrix} - \rank\mathcal{F}_{n,
[\mathscr{E},\mathscr{A}]}.
\end{eqnarray*}
Thus, it follows from Proposition \ref{prop:ker} that rank condition \eqref{thm:partialcausalityb} holds if, and only if, 
\begin{eqnarray}\label{eq:causal2}
\ker{\begin{bmatrix} \mathcal{Z}\mathcal{A} \\ \mathcal{C} \\  \mathcal{F}_{n,[\bar{E},\bar{A}]} \end{bmatrix}} \subseteq \ker{\mathcal{K}} .
\end{eqnarray}
Since $\ker(\mathcal{Z}\mathcal{A})  = {\mathcal{A}}^{-1} \left(\ker \mathcal{Z}\right) = {\mathcal{A}}^{-1} \left(\im{\mathcal{F}_{n, [\mathscr{E},\mathscr{A}]}}\right)$, \eqref{eq:causal2} is equivalent to
\begin{equation}\label{eq:sub1}
\mathcal{A}^{-1}\left(\im{\mathcal{F}_{n,[\mathscr{E},\mathscr{A}]}}\right) \cap \ker{\mathcal{C}} \cap \ker \mathcal{F}_{n,[\bar{E},\bar{A}]}  \subseteq \ker{\mathcal{K}}.
\end{equation}

$(ii) \Rightarrow (iii)$: Let $v_n\in \mathcal{A}_1^{-1}(\im{\mathcal{F}_{n,[\mathscr{E},\mathscr{A}]}})\cap \mathcal{W}^{*}_{[E,A,0,C]}$ be arbitrary. By~\eqref{eq:W*-barEA-EAC} we find that $\mathcal{W}^{*}_{[E,A,0,C]} = \mathcal{W}^*_{[\bar{E}, \bar{A},0,0]} \cap \ker C$ and since the Wong sequences terminate after finitely many steps and in each iteration before termination the dimension increases by at least one, we have $\mathcal{W}^*_{[\bar{E}, \bar{A},0,0]} = \mathcal{W}^{n}_{[\bar{E}, \bar{A},0,0]}$. Therefore, $v_n\in \mathcal{A}_1^{-1}(\im{\mathcal{F}_{n,[\mathscr{E},\mathscr{A}]}})\cap \ker C \cap \mathcal{W}^{n}_{[\bar{E}, \bar{A},0,0]}$. Hence, in particular, there exist $v_i \in \mathcal{W}^{i}_{[\bar{E}, \bar{A},0,0]}$,  $i=1,\ldots,n-1$, such that
\[
 \bar{E}v_{n} + \bar{A}v_{n-1} = 0,~ \bar{E}v_{n-1} + \bar{A}v_{n-2} = 0,~ \ldots,~ \bar{E}v_{2} + \bar{A}v_{1} = 0,~\bar{E}v_1 = 0.
\]
This implies that $\mathcal{F}_{n,[E_{\epsilon}, A_{\epsilon}]}v = 0$ for $v=(v_n^\top,\ldots,v_1^\top)^\top$. Furthermore, we have that
\begin{subequations}\label{eq:sub7}
	\begin{align}
\mathcal{A}^{-1}(\im{\mathcal{F}_{n,[\mathscr{E},\mathscr{A}]}}) &= \mathcal{A}_1^{-1}(\im{\mathcal{F}_{n,[\mathscr{E},\mathscr{A}]}}) \times \overset{(n-1) \text{-times}}{\overbrace{\mathbb{R}^n \times \ldots \times \mathbb{R}^n}}, \\ \ker{\mathcal{K}} &= \ker
\NiceMatrixOptions
{nullify-dots,code-for-last-col = \color{black},code-for-last-col=\color{black}}
\begin{bNiceMatrix}[first-row]
& \Ldots[line-style={solid,<->},shorten=0pt]^{n \text{-block columns}} \\
K & 0 & \ldots & 0
\end{bNiceMatrix}
=
\ker{K} \times \overset{(n-1) \text{-times}}{\overbrace{\mathbb{R}^n \times \ldots \times \mathbb{R}^n}}, \\
\ker{\mathcal{C}} &= \ker
\NiceMatrixOptions
{nullify-dots,code-for-last-col = \color{black},code-for-last-col=\color{black}}
\begin{bNiceMatrix}[first-row]
& \Ldots[line-style={solid,<->},shorten=0pt]^{n \text{-block columns}} \\
C & 0 & \ldots & 0
\end{bNiceMatrix}
=
\ker{C} \times \overset{(n-1) \text{-times}}{\overbrace{\mathbb{R}^n \times \ldots \times \mathbb{R}^n}},
\end{align}
\end{subequations}
from which it follows
\[
    v \in \mathcal{A}^{-1}\left(\im{\mathcal{F}_{n,[\mathscr{E},\mathscr{A}]}}\right) \cap \ker{\mathcal{C}} \cap \ker \mathcal{F}_{n,[\bar{E},\bar{A}]}  \subseteq \ker{\mathcal{K}},
\]
hence $v_n \in \ker K$.

$(ii) \Leftarrow (iii)$: If $v \in \mathcal{A}^{-1}\left(\im{\mathcal{F}_{n,[\mathscr{E},\mathscr{A}]}}\right) \cap \ker{\mathcal{C}} \cap \ker \mathcal{F}_{n,[\bar{E},\bar{A}]}$, then, with a similar argument as in the previous step, for $v=(v_n^\top,\ldots,v_1^\top)^\top$ it follows that $v_i \in \mathcal{W}^{i}_{[\bar{E}, \bar{A},0,0]}$,  $i=1,\ldots,n$; in particular $v_n \in \mathcal{W}^n_{[\bar{E}, \bar{A},0,0]} = \mathcal{W}^*_{[\bar{E}, \bar{A},0,0]}$. Then invoking \eqref{eq:sub7} and \eqref{eq:W*-barEA-EAC} it follows that $v_n\in \mathcal{A}_1^{-1}(\im{\mathcal{F}_{n,[\mathscr{E},\mathscr{A}]}})\cap \mathcal{W}^{*}_{[E,A,0,C]} \subseteq \ker K$, thus $v\in \ker{\mathcal{K}}$.

$(iii) \Leftrightarrow (iv)$:
In order to prove this, it is sufficient to show
\begin{equation}\label{eq:pf8}
A^{-1}\left(E \left( \mathcal{V}_{[E,A,B,0]}^{n-1} \right) \right)
  = \mathcal{A}_1^{-1}(\im{ \mathcal{F}_{n,[\mathscr{E},\mathscr{A}]}}). 
\end{equation}
Let a nonzero vector $z \in A^{-1}\left(E \left( \mathcal{V}_{[E,A,B,0]}^{n-1} \right) \right)$
be given.
Then, there exists $v_{n-1} \in \mathcal{V}_{[E,A,B,0]}^{n-1}$ 
such that $Az = -Ev_{n-1}$.
Therefore, there exist $v_{i} \in \mathcal{V}^{i}_{[E,A,B,0]}$ and $u_{i+1} \in \mathbb{R}^l$, for $0 \leq i \leq n-2$, such that
\begin{subequations}\label{eq:pf1}
	\begin{eqnarray}
Ev_{i-1} +Av_{i} + Bu_{i} &=& 0 , ~\text{ for } 1 \leq i \leq n-1\\
Ev_{n-1} + Az &=& 0 .	
	\end{eqnarray}
\end{subequations}
By taking $v = \left[\begin{smallmatrix} \bar{v}_0 \\ \vdots \\ \bar{v}_{n-1} \end{smallmatrix}\right]$, $\bar{v}_i = \begin{bmatrix} v_i \\ u_i \end{bmatrix}$, for $0 \leq i \leq n-1$, where $u_0:=0$,
and using the definitions of   $\mathcal{F}_{n,[\mathscr{E},\mathscr{A}]}$ and $\mathcal{A}_1$,
system \eqref{eq:pf1} can be rewritten as,
\begin{eqnarray*}\label{eq:pf4}
\mathcal{F}_{n,[\mathscr{E},\mathscr{A}]}v + \mathcal{A}_1z = 0,\quad  \emph{i.e.},\quad  z \in \mathcal{A}_1^{-1}(\im{\mathcal{F}_{n,[\mathscr{E},\mathscr{A}]}}).
\end{eqnarray*}
Thus,
\begin{equation}\label{eq:pf6}
A^{-1}\left(E \left( \mathcal{V}_{[E,A,B,0]}^{n-1} \right) \right) \subseteq \mathcal{A}_1^{-1}(\im{\mathcal{F}_{n,[\mathscr{E},\mathscr{A}]}}).
\end{equation}

Now, let a nonzero vector $z \in \mathcal{A}_1^{-1}(\im{\mathcal{F}_{n,[\mathscr{E},\mathscr{A}]}})$ 
be given. This implies that for some vector $v = \left[\begin{smallmatrix} \bar{v}_0 \\ \vdots \\ \bar{v}_{n-1} \end{smallmatrix}\right]$, where $\bar{v}_i = \begin{bmatrix}
	v_i \\ u_i
\end{bmatrix} \in \mathbb{R}^{n+l}$ for $i \in \{0,1,\ldots,n-1 \}$, we have $\mathcal{A}_1z = -\mathcal{F}_{n,[\mathscr{E},\mathscr{A}]}v$.
Using the definitions of   $\mathcal{F}_{n, [\mathscr{E},\mathscr{A}]}$ and $\mathcal{A}_1$, the system $\mathcal{A}_1z + \mathcal{F}_{n,[\mathscr{E},\mathscr{A}]}v = 0$ can be written as \eqref{eq:pf1}. Therefore, it follows from the definition of the sequence $\left\{\mathcal{V}^{i}_{[E,A,B,0]} \right\}_{i = 0}^{\infty}$ 
that $v_i \in \mathcal{V}^{i}_{[E,A,B,0]} $, for $i \in \{0,1,\ldots,n-1\}$, and $Ev_{n-1} = Az$.
Therefore,
$z \in A^{-1}\left(E \left( \mathcal{V}_{[E,A,B,0]}^{n-1} \right) \right)$, 
and hence
\begin{equation}\label{eq:pf7}
\mathcal{A}_1^{-1}(\im{\mathcal{F}_{n,[\mathscr{E},\mathscr{A}]}}) \subseteq 
A^{-1}\left(E \left( \mathcal{V}_{[E,A,B,0]}^{n-1} \right) \right).
\end{equation}
Thus, \eqref{eq:pf8} follows from \eqref{eq:pf6} and \eqref{eq:pf7}.

$(i) \Leftrightarrow (v)$:
In view of the rank identities  \eqref{eq:barFn} and \eqref{eq:barFn1}, \eqref{thm:partialcausalityb} is equivalent to \eqref{eq:Fbar}. Now, it follows from Step $3$ of the proof of Theorem \ref{thm:partial:causal:1} that rank condition \eqref{eq:Fbar} holds if, and only if,
\begin{equation}\label{eq:Kcausal_O}
\rank \mathcal{F}_{n+1,[\bar{E}_O,\bar{A}_O]} = \rank \mathcal{F}_{n+1,[\bar{E}_O, \bar{A}_O, K_O]}.
\end{equation}
Here, $\bar{E}_O = \begin{bmatrix}
E_O \\ 0 \end{bmatrix}$, $\bar{A}_O = \begin{bmatrix} A_O \\ C_O \end{bmatrix}$, $\bar{B}_O = \begin{bmatrix} B_O & 0 \\ D & -I_p \end{bmatrix}$, $E_O$, $A_O$, $B_O$ correspond to the decomposition \eqref{Cdecomposition} of $E$, $A$, $B$, and $C_O$, $K_O$ are the corresponding parts of $C$, $K$ according to the decomposition \eqref{Cdecomposition}, respectively. In addition, by Proposition \ref{prop:kcd}, for the tuple $\{ E_O,A_O,B_O,C_O,D\}$ there exist two nonsingular matrices $\tilde{U}$ and $\tilde{V}$ such that
\begin{equation}\label{eq:completelycont}
\tilde{U}E_O\tilde{V} = \begin{bmatrix}
E_{11} & E_{12} \\ 0 & I_{m_2}
\end{bmatrix}, ~\tilde{U}A_O\tilde{V} = \begin{bmatrix} A_{11} & A_{12} \\ 0 & A_{22} \end{bmatrix},
~\tilde{U}B_O = \begin{bmatrix}
B_1 \\ 0 \end{bmatrix},~ C_O\tilde{V} = \begin{bmatrix} C_{11} & C_{12} \end{bmatrix},
\end{equation}
where  $\{E_{11},A_{11},B_{1},C_{11},D\}$ represents the completely controllable part of~\eqref{dls} and $m_2 \in \mathbb{N} \cup \{0\}$.

Thus, to prove the equivalence of statements $(i)$ and $(v)$, it is sufficient to show that condition \eqref{eq:Kcausal_O} is equivalent to partial impulse observability of  $\{E_{11},A_{11},C_{11}\}$ with respect to $K_{11}$, where $K_O\tilde{V} = \begin{bmatrix}
	K_{11} & K_{12} & K_{13}
\end{bmatrix}$. Now, set
\begin{eqnarray*}
&\mathcal{U}_1 = \begin{bNiceMatrix}
\blkdiag\{\tilde{U},I_p\} &  &  \\
& \Ddots^{(n+1) \text{ times}} & \\
& & \blkdiag\{\tilde{U},I_p\}
\end{bNiceMatrix} ,~\mathcal{V}_1 =   \begin{bNiceMatrix}
V &  &  \\
& \Ddots^{(n+1) \text{ times}} & \\
& & V  \end{bNiceMatrix} ,~\mathcal{U}_2 = \begin{bmatrix}
\mathcal{U}_1 \\ & I_r
\end{bmatrix} .&
\end{eqnarray*}
Clearly, $\rank \mathcal{F}_{n+1,[\bar{E}_O,\bar{A}_O]} = \rank (\mathcal{U}_1 \mathcal{F}_{n+1,[\bar{E}_O,\bar{A}_O]} \mathcal{V}_1)$. We now write the matrix $\mathcal{F}_{n+1,[\bar{E}_O,\bar{A}_O]}$ in terms of $E_O$, $A_O$, $C_O$, and obtain all the $2(n+1)$-block rows of the matrix $\mathcal{U}_1 \mathcal{F}_{n+1,[\bar{E}_O,\bar{A}_O]} \mathcal{V}_1$. Thus, substituting decomposition \eqref{eq:completelycont} in all block rows of $\mathcal{U}_1 \mathcal{F}_{n+1,[\bar{E}_O,\bar{A}_O]} \mathcal{V}_1$, we see that an identity matrix $I_{m_2}$ appears $(n+1)$-times on the diagonal. By permuting those matrices to the upper left corner and applying Proposition \ref{prop3}, we obtain
\begin{equation*}
\rank \mathcal{F}_{n+1,[\bar{E}_O,\bar{A}_O]} = (n+1) \rank{I_{m_2}} + \rank \mathcal{F}_{n+1,[\bar{E}_{11},\bar{A}_{11}]},
\end{equation*}
where $\bar{E}_{11} = \begin{bmatrix}
E_{11} \\ 0 \end{bmatrix}$ and $\bar{A}_{11} = \begin{bmatrix}
A_{11} \\ C_{11} \end{bmatrix}$. In a similar manner, we obtain
\begin{eqnarray*}
\rank \mathcal{F}_{n+1,[\bar{E}_O,\bar{A}_O,K_O]} &=& \rank \left( \mathcal{U}_2 \mathcal{F}_{n+1,[\bar{E}_O,\bar{A}_O,K_O]} \mathcal{V}_1 \right) \\
&=& (n+1) \rank{I_{m_2}} + \rank \mathcal{F}_{n+1,[\bar{E}_{11},\bar{A}_{11},K_{11}]}.
\end{eqnarray*}
Thus, it follows from Proposition \ref{prop:ker} that rank identity \eqref{eq:Kcausal_O} is equivalent to
\begin{eqnarray}\label{eq:causal3}
\ker \mathcal{F}_{n+1,[\bar{E}_{11},\bar{A}_{11}]} \subseteq \ker
\NiceMatrixOptions
{nullify-dots,code-for-last-col = \color{black},code-for-last-col=\color{black}}
\begin{bNiceMatrix}[first-row]
& \Ldots[line-style={solid,<->},shorten=0pt]^{(n+1) \text{-block columns}} \\
0 & K_{11} & 0 & \cdots & 0 \end{bNiceMatrix}. 
\end{eqnarray}
We show that \eqref{eq:causal3} is equivalent to 
\begin{equation}\label{eq:causal4}
\mathcal{W}^*_{[E_{11},A_{11},0,C_{11}]} \cap A_{11}^{-1}(\im{E_{11}})
\subseteq \ker K_{11}.
\end{equation} 
To see ``$\Leftarrow$'', let $v = (v_{n+1}^\top,\ldots,v_1^\top)^\top \in \ker \mathcal{F}_{n+1,[\bar{E}_{11},\bar{A}_{11}]}$, then
\begin{equation}\label{eq:kerFn-barE11A11}
\bar{E}_{11}v_1 = 0,~ \bar{E}_{11}v_{i+1} + \bar{A}_{11}v_{i} = 0, \text{ for } 1 \leq i \leq n.
\end{equation}
In particular, $v_n \in \mathcal{W}^{n}_{[\bar{E}_{11}, \bar{A}_{11},0,0]}$ and since $\bar{E}_{11}v_{n+1} + \bar{A}_{11}v_{n} = 0$ we further have $v_n \in \ker C_{11} \cap A_{11}^{-1}(\im{E_{11}})$. Again, since the Wong sequences terminate after finitely many steps and the dimension increases in each step, we have $\mathcal{W}^{n}_{[\bar{E}_{11}, \bar{A}_{11},0,0]} = \mathcal{W}^{*}_{[\bar{E}_{11}, \bar{A}_{11},0,0]}$, and from~\eqref{eq:W*-barEA-EAC} it follows that $v_n \in \mathcal{W}^*_{[E_{11},A_{11},0,C_{11}]} \cap A_{11}^{-1}(\im{E_{11}})
\subseteq \ker K_{11}$, thus $v \in \ker \begin{bmatrix} 0 & K_{11} & 0 & \cdots & 0\end{bmatrix}$.

For ``$\Rightarrow$'', let $v_n \in \mathcal{W}^*_{[E_{11},A_{11},0,C_{11}]} \cap A_{11}^{-1}(\im{E_{11}})$. Then, with a similar argument as in the previous step,  $v_n \in \mathcal{W}^{n}_{[\bar{E}_{11}, \bar{A}_{11},0,0]} \cap \ker C_{11} \cap A_{11}^{-1}(\im{E_{11}})$, hence there exist $v_{n+1}\in\mathbb{R}^{n_1}$ and $v_i\in \mathcal{W}^{i}_{[\bar{E}_{11}, \bar{A}_{11},0,0]}$, $i=1,\ldots,n$, such that~\eqref{eq:kerFn-barE11A11} holds, thus $v = (v_{n+1}^\top,\ldots,v_1^\top)^\top \in \ker \mathcal{F}_{n+1,[\bar{E}_{11},\bar{A}_{11}]}\subseteq \ker \begin{bmatrix} 0 & K_{11} & 0 & \cdots & 0\end{bmatrix}$, by which $v_n\in \ker K_{11}$.

Notably, \eqref{eq:causal4} is equivalent to partial impulse observability of $\{E_{11},A_{11},C_{11}\}$ with respect to $K_{11}$, cf. Proposition \ref{prop:partial:impulse}. This completes the proof.
\end{proof}

In view of the above results, the following remark is warranted.
\begin{remark}\label{rem:K=I}
For $K = I_n$, the statement $(iii)$ in Theorem \ref{thm:Kcausallydetectable}  reduces to
\begin{equation}\label{eq:sub3}
\mathcal{A}_1^{-1}(\im{\mathcal{F}_{n,[\mathscr{E},\mathscr{A}]}}) \cap \mathcal{W}^{*}_{[E,A,0,C]} = \{0\}.
\end{equation}
Since $\mathcal{W}^{*}_{[E,A,0,C]} = \bigcup_{i \in \mathbb{N}} \mathcal{W}^{i}_{[E,A,0,C]}$, \eqref{eq:sub3} implies
\begin{eqnarray}\label{eq:sub4} \mathcal{A}_1^{-1}(\im{ \mathcal{F}_{n,[\mathscr{E},\mathscr{A}]}}) \cap \mathcal{W}^{1}_{[E,A,0,C]} = \{0\}.
\end{eqnarray}
Further, by definition of the generalized Wong sequences, $\mathcal{W}^{1}_{[E,A,0,C]} = \ker{E} \cap \ker{C}$. Therefore, \eqref{eq:sub4} becomes
$ \mathcal{A}_1^{-1}(\im{ \mathcal{F}_{n,[\mathscr{E},\mathscr{A}]}}) \cap \ker{C} \cap \ker{E} = \{0\}$.
Thus in this case, Theorem \ref{thm:partialcausality} implies {\it causal detectability} of system \eqref{dls}, which is necessary and sufficient for the full-state estimation via system \eqref{obsv}; for more details, see \cite[Thm.~1]{jaiswal2021necessary}. Likewise, again invoking $\ker{E} \cap \ker{C} \subseteq \mathcal{W}^{*}_{[E,A,0,C]}$, the characterizations $(iv)$ and $(v)$ in Theorem \ref{thm:Kcausallydetectable}, for the case $K=I$, imply alternative characterizations for {\it causality} of system \eqref{dls}, which can be found in \cite{berger2019ode,jaiswal2021necessary}.
\end{remark}

\section{Functional ODE estimator}\label{sec:ODE:estimator}
In this section, we will prove that partial causal detectability of system \eqref{dls} with respect to $K$ is necessary and sufficient for the estimation of the functional vector $z(t)$ in \eqref{dls} via system \eqref{obsv}. First, we exploit the behavior $\mathscr{B}$ to give a precise definition of functional ODE estimators for \eqref{dls}, similar to \cite[Def.~3.2]{berger2019disturbance}.
\begin{definition}\label{def:observer}
System \eqref{obsv} is said to be a functional ODE estimator for \eqref{dls}, if for every  $(x,u,y,z) \in \mathscr{B}$ there exist $w\in \mathcal{AC}_{\loc}(\mathbb{R};\mathbb{R}^l)$ and $\hat z\in \mathscr{L}^1_{\loc}(\mathbb{R};\mathbb{R}^r)$ such that $(w,u,y,\hat z)$ satisfy \eqref{obsv} for almost all $t\in\mathbb{R}$, and for all $w,\hat z$ with this property, $ \hat{z}(t) \to z(t)$ for $t\to\infty$.
\end{definition}

\begin{remark}\label{rem:obsv}
Note that if a functional ODE estimator satisfies the state matching property, \emph{i.e.}, $\hat{z}(0) = z(0)$ implies $\hat{z}(t) = z(t)$, for almost all $t > 0$, then it is known as a functional ODE observer. In case $K=I_n$, this condition holds automatically and, therefore, there is no difference between ODE observer and ODE estimator. However, in the case of partial-state estimation (\emph{i.e.}, $K \neq I_n$), the state matching condition is not always necessary to hold by default. Therefore, ODE observer and ODE estimator are not the same in case of partial-state estimation. We will show this fact in Example \ref{exp2} below.
\end{remark}

Before providing the main result of this section, we will establish a necessary condition for partial-state estimation of the $\sigma$-block in the QKF \eqref{eq:qkf} of \eqref{dls} by a functional ODE estimator \eqref{obsv}.
\begin{lemma}\label{lm:sigma}
Consider the system \begin{subequations}\label{dls:1}
\begin{eqnarray}
J_{\sigma}\dot{x}_{\sigma}(t) &=& x_{\sigma}(t)+B_{\sigma}u(t), \label{dls:1a} \\
y_{\sigma}(t) &=& 0, \label{dls:1b}\\
z_{\sigma}(t) &=& K_{\sigma} x_{\sigma}(t), \label{dls:1c}
\end{eqnarray}
\end{subequations}
where $J_{\sigma}$ is a nilpotent matrix with nilpotency index $h$. If there exists a functional ODE estimator \eqref{obsv} for system \eqref{dls:1}, then $K_{\sigma}J_{\sigma}^{i}B_{\sigma} = 0$ for all $1 \leq i \leq h$.
\end{lemma}

\begin{proof}
Assume that there exist a functional ODE estimator for the system \eqref{dls:1}. Then the estimator is given by
\begin{subequations}\label{obsv1}
\begin{eqnarray}
\dot{w}(t) &=& Nw(t) + H u(t), \\
\hat{z}_{\sigma}(t) &=& Rw(t) +  M u(t),
\end{eqnarray}
\end{subequations}
and, by S$2$), the estimate $\hat z_{\sigma}$ is given by
\begin{equation*}
\hat{z}_{\sigma}(t) = R\left( \exp(Nt)w(0) + \int_0^t \exp(N(t-\tau))Hu(\tau) \rm{d} \tau \right) + Mu(t).
\end{equation*}
Also, by S$3$), the solution of the system \eqref{dls:1} is given by
\begin{equation*}
z_{\sigma}(t) = - \sum_{i = 0}^{h - 1} K_{\sigma}J_{\sigma}^i B_{\sigma} u^{(i)}(t).
\end{equation*}
Since system \eqref{obsv1} is a functional ODE estimator for system \eqref{dls:1}, we have $e(t) := \hat{z}_{\sigma}(t) - z_{\sigma}(t) \to 0$ as $t \to \infty$ for each input function $u$ and initial value $w(0)$.

Let $s$ be the largest index  such that $K_{\sigma} J_\sigma^s B_{\sigma} \neq 0$ for $1 \leq s \leq h-1$. Choose $w(0) = 0$ and $u(t) = \frac{\sin(t^2)}{t^s}  e_k$ with $e_k$ being an arbitrary unit vector for $1 \leq k \leq m$. Then it is straightforward to see that
$u^{(i)}(t) \to 0 \text{ for } i=0,\ldots,s-1 \text{ and } u^{(s)}(t) \not\to 0$.
Since $R \exp(Nt) \to 0$ (which can be seen from choosing $u=0$ and arbitrary $w(0)$), it is easy to show that $\int_0^t R \exp(N(t-\tau)) u(\tau) d\tau \to 0$ and together with $e(t)\to 0$ it follows that $K_\sigma J_\sigma^s B_\sigma u^{(s)}(t) \to 0$, which is only possible when  $K_\sigma J_\sigma^s B_\sigma e_k = 0$. Since~$k$ was arbitrary it follows that $K_\sigma J_\sigma^s B_\sigma = 0$, which contradicts the assumption on the index $s$. Therefore, $K_\sigma J_\sigma^i B_\sigma = 0$ for all $i=1,\ldots,h-1$. This completes the proof.
\end{proof}

In the following theorem, we prove that partial causal detectability of system~\eqref{dls} is equivalent to the existence of a functional ODE estimator.
\begin{theorem}\label{thm:ODE:estimator}
For system \eqref{dls}, the following statements are equivalent:
\begin{enumerate}
\item[(i)] System \eqref{dls} is partially causal detectable with respect to $K$. \label{bi}
		
\item[(ii)] There exists a functional ODE estimator for system \eqref{dls}. \label{ai}
\end{enumerate}
\end{theorem}

\begin{proof}
$(i) \Rightarrow (ii)$: To prove this part, first, we give a step-by-step procedure to design a functional ODE estimator of the form \eqref{obsv}. \\
\textbf{Step $1$:} Compute orthogonal matrices $U_O$ and $V_O$ according to Lemma \ref{lm:Cdecomposition}, which transform $\{E,A,B\}$ as in \eqref{Cdecomposition}, and obtain $\{ E_O,A_O,B_O \}$. Define $$CV_O = \begin{bmatrix}
	C_O & C_{k-1} & \hdots & C_1
\end{bmatrix} \text{ and } KV_O = \begin{bmatrix}
	K_O & K_{k-1} & \hdots & K_1
\end{bmatrix}.$$
\textbf{Step $2$:} According to Lemma \ref{lm:qkf}, compute nonsingular matrices $P$ and $Q$ such that $(\lambda \bar{E}_O -\bar{A}_O)$ is in QKF \eqref{eq:qkf}, \emph{i.e.},
\begin{eqnarray*}
&P(\lambda\bar{E}_O -\bar{A}_O)Q = \blkdiag\{ \lambda E_{\epsilon}- A_{\epsilon}, \lambda I_{n_f} - J_f , \lambda J_{\sigma} - I_{n_\sigma} ,\lambda E_{\eta} - A_{\eta} \},~& \\
& P\bar{B}_O := \begin{bmatrix}
	B_{\epsilon}^\top & B_f^\top & B_{\sigma}^\top & B_{\eta}^\top \end{bmatrix}^\top,\text{ and} ~~K_OQ :=\begin{bmatrix}
K_{\epsilon} & K_f & K_{\sigma} & K_{\eta}
\end{bmatrix},
\end{eqnarray*}
where $\bar{E}_O = \begin{bmatrix}
E_O \\ 0 \end{bmatrix}$, $\bar{A}_O = \begin{bmatrix} A_O \\ C_O \end{bmatrix}$, and $\bar{B}_O = \begin{bmatrix} B_O & 0 \\ D & -I_p \end{bmatrix}$. \\
\textbf{Step $3$:} Utilizing the Jordan decomposition, compute a non-singular matrix $U_1$ such that
$U_1^{-1} J_fU_1 = \blkdiag\{
		J_{f_1}, J_{f_2}\}$,
where $\sigma(J_{f_1}) \subseteq \overline{\mathbb{C}^+}$ and $\sigma(J_{f_2}) \subseteq \mathbb{C}^-$. Set $U_1^{-1}B_f = \begin{bmatrix}
			B_{f_1} \\ B_{f_2}
\end{bmatrix}$ and $K_fU_1 = \begin{bmatrix}
			K_{f_1} & K_{f_2}
\end{bmatrix}$. \\
\textbf{Step $4$:} Utilizing the singular value decomposition, compute a nonsingular matrix~$U_2$ such that $U_2E_{\eta} = \begin{bmatrix}
I_{n_\eta} \\ 0 \end{bmatrix}$.
Set $U_2A_{\eta} = \begin{bmatrix} A_{\eta_1} \\ A_{\eta_2} \end{bmatrix}$ and $U_2B_{\eta} = \begin{bmatrix} 
B_{\eta_1} \\ B_{\eta_2} \end{bmatrix}$. \\
\textbf{Step $5$:} Set $x=V_O\begin{bmatrix}
x_k^\top & x_{k-1}^\top & \ldots & x_1^\top
\end{bmatrix}^\top$, $\bar{u} := \begin{bmatrix}	
u \\ y
\end{bmatrix}$, and
\begin{equation*}
x_k := \blkdiag\{I_{n_\epsilon},U_f,I_{n_\sigma},I_{n_\eta} \} Q\begin{bmatrix}
x_\epsilon^\top & x_{f_1}^\top & x_{f_2}^\top & x_{\sigma}^\top & x_\eta^\top
\end{bmatrix}^\top .
\end{equation*} 
In the new coordinates, system \eqref{dls} becomes
\begin{eqnarray*}
E_{\epsilon} \dot{x}_{\epsilon}(t) &=& A_{\epsilon}x_{\epsilon}(t) + B_{\epsilon} \bar{u}(t), \label{dls3a} \\
\dot{x}_{f_1}(t) &=& J_{f_1}x_{f_1}(t) + B_{f_1}\bar{u}(t), \label{dls3b} \\
\dot{x}_{f_2}(t) &=& J_{f_2}x_{f_2}(t) + B_{f_2}\bar{u}(t), \label{dls3c} \\
J_{\sigma} \dot{x}_{\sigma}(t) &=& x_{\sigma}(t) + B_{\sigma}\bar{u}(t), \label{dls3d} \\
\dot{x}_{\eta}(t) &=& A_{\eta_1} x_{\eta}(t) + B_{\eta_1}\bar{u}(t), \label{dls3e} \\
0 &=& A_{\eta_2} x_{\eta}(t) + B_{\eta_2}\bar{u}(t), \label{dls3f} \\
z(t) &=& K_{\epsilon}x_{\epsilon}(t) + K_{f_1}x_{f_1}(t) + K_{f_2}x_{f_2}(t) + K_{\sigma}x_{\sigma}(t) + K_{\eta}x_{\eta}(t). \label{dls3g}
\end{eqnarray*}
Here $x_1 = x_2 = \ldots = x_{k-1} = 0$ due to decomposition \eqref{Cdecomposition}, for details see Step $2$ in the proof of Theorem \ref{thm:partial:causal:1}. \\
\textbf{Step $6$:} As shown in Step $2$ of the proof of Theorem \ref{thm:partial:causal:1}, partial detectability of $\{E,A,C\}$ with respect to $K$ implies that $\{E_O,A_O,C_O\}$ is partially detectable with respect to $K_O$.  Hence it follows from \cite[Lem. $4$]{jaiswal2023detectability} that $K_{\epsilon} = 0$ and $K_{f_1} = 0$. \\
\textbf{Step $7$:} The solution of the $\sigma$-block is given by $x_{\sigma}(t) = - \sum_{i = 0}^{h}J_{\sigma}^i B_{\sigma}\bar{u}^{(i)}(t) $ and the tuple $\{\bar{E}_O,\bar{A}_O,\bar{B}_O\}$ is partially causal with respect to $K_O$, since $\{\bar E, \bar A, \bar B\}$ is partially causal with respect to $K$ by assumption. So, \eqref{eq:Sigma} reveals that $K_{\sigma}J_{\sigma} = 0$ and hence, $K_{\sigma}x_{\sigma}(t) = - K_{\sigma}B_{\sigma}\bar{u}(t)$. \\	
\textbf{Step $8$:} In the new coordinates, the problem of functional ODE estimator design for system \eqref{dls} reduces to the problem of functional ODE estimator design for
\begin{subequations}\label{dls5}
\begin{eqnarray}
\dot{x}_{f_2}(t) &=& J_{f_2}x_{f_2}(t) + B_{f_2}\bar{u}(t), \label{dls5a} \\
\dot{x}_{\eta}(t) &=& A_{\eta_1} x_{\eta}(t) + B_{\eta_1}\bar{u}(t), \label{dls5b} \\
0 &=& A_{\eta_2} x_{\eta}(t) + B_{\eta_2}\bar{u}(t), \label{dls5c} \\
z(t) &=& K_{f_2}x_{f_2}(t) + K_{\eta}x_{\eta}(t) - K_{\sigma}B_{\sigma} \bar{u}(t). \label{dls5d}
\end{eqnarray}
\end{subequations}
\textbf{Step $9$:} Since $\rank \begin{bmatrix} \lambda I_{n_\eta} - A_{\eta_1}  \\ -A_{\eta_2} \end{bmatrix}=n_\eta$ for all $\lambda \in\mathbb{C}$ by Lemma \ref{lm:qkf}, there exists $L\in\mathbb{R}^{n_\eta \times (m_\eta - n_\eta)}$ such that $\sigma(A_{\eta_1} - L A_{\eta_2}) \subseteq \mathbb{C}^-$. \\
\textbf{Step $10$:} We claim that the following system is a functional ODE estimator for \eqref{dls5}:
\begin{eqnarray*}
\dot{w}(t) &=& Nw(t) + H\bar{u}(t), \\ 
\hat{z}(t) &=& Rw(t) + M\bar{u}(t), 
\end{eqnarray*}
where $N = \blkdiag\{J_{f_2}, A_{\eta_1} - L A_{\eta_2} \}$,  $R = \begin{bmatrix} K_{f_2} & K_{\eta}	\end{bmatrix}$, $M = -K_{\sigma}B_{\sigma}$, and $H = \begin{bmatrix} B_{f_2} \\ B_{\eta_1} - L B_{\eta_2} \end{bmatrix}$. Set $e := \hat{z} - z$ and $e_1 := w - \begin{bmatrix}
 x_{f_2} \\ x_{\eta}
\end{bmatrix}$. Then
\begin{eqnarray*}
\dot{e}_1(t) &=& Ne_1(t) + \begin{bmatrix}
	0 \\ L(A_{\eta_2}x_{\eta}(t) + B_{\eta_2}\bar{u} (t))
\end{bmatrix} = Ne_1(t), \\
e(t) &=& Re_1(t).
\end{eqnarray*}
Since $\sigma(N) \subseteq \mathbb{C}^-$, $e_1(t) \to 0$ as $t \to \infty$. Consequently, $e(t) \to 0$ as $t \to \infty$.

$(ii) \Rightarrow (i)$: Assume that system \eqref{dls} has a functional ODE estimator. Then, with the same proof as in \cite[Thm. 2]{jaiswal2023detectability}, partial detectability with respect to $K$ can be inferred. Now, by repeating Step $1$ to Step $6$ of the first part of the proof, we obtain the system in the following form
\begin{eqnarray*}
\dot{x}_{f_2}(t) &=& J_{f_2}x_{f_2}(t) + B_{f_2}\bar{u}(t), \\
J_{\sigma} \dot{x}_{\sigma}(t) &=& x_{\sigma}(t) + B_{\sigma}\bar{u}(t), \\
\dot{x}_{\eta}(t) &=& A_{\eta_1} x_{\eta}(t) + B_{\eta_1}\bar{u}(t),  \\
0 &=& A_{\eta_2} x_{\eta}(t) + B_{\eta_2}\bar{u}(t),  \\
z(t) &=& K_{f_2}x_{f_2}(t) + K_{\eta}x_{\eta}(t) + K_{\sigma}x_{\sigma}(t).
\end{eqnarray*}
By the definition of functional ODE estimators, if one exists for the above system, then also one exists for the system \eqref{dls:1}. Hence, it follows from Lemma \ref{lm:sigma} that $K_{\sigma}J_{\sigma}^iB_{\sigma} = 0$ for all $i \geq 1$. Since the QKF \eqref{eq:qkf} is computed for the triple $\{\bar{E}_O,\bar{A}_O,\bar{B}_O\}$ and $\begin{bmatrix}
\bar{E}_O & \bar{B}_O
\end{bmatrix}$ has full row rank (see Step $3$ in the first part), by performing a similar calculation as done in the proof of Theorem \ref{thm:partial:causal:1}, it is easy to conclude that $\begin{bmatrix}
J_{\sigma} & B_{\sigma}
\end{bmatrix}$
has also full row rank.
By repeating the same steps as done in Step $5$ in the proof of Theorem \ref{thm:partial:causal:1}, we obtain that $K_{\sigma}J_{\sigma} = 0$. Thus, by Definition \ref{def:partial:causality} and the solution discussion in S$1$)-S$4$),  $\{\bar{E}_O,
\bar{A}_O,\bar{B}_O\}$ is partially causal with respect to $K_O$. Therefore, Step $2$ in the proof of Theorem \ref{thm:partial:causal:1} implies that $\{\bar{E},\bar{A}, \bar{B}\}$ is partially causal with respect to $K$. This completes the proof.
\end{proof}

\section{Numerical illustration}\label{sec:numerical}
In this section a numerical example is given to illustrate the theoretical findings.  Also, Example \ref{exp2} reveals that it is not always possible to design a functional ODE observer, if a functional ODE estimator exists for the system \eqref{dls}.

\begin{example}\label{exp2}
Consider system \eqref{dls} with coefficient matrices \\
$  E = \begin{bmatrix}
1 & 0 & 0 & 0 \\ 0 & 1 & 0 & 0 \\
0 & 0 & 1 & 0 \\ 0 & 0 & 0 & 0
\end{bmatrix},~A = \begin{bmatrix}
	1 & 0 & 0 & 0 \\ 0 & -1 & 1 & 0 \\
	0 & 0 & -1 & 0 \\ 0 & 0 & 0 & 1
\end{bmatrix},~B = \begin{bmatrix}
	1 \\ 1 \\ 1 \\ 1
\end{bmatrix},~C = \begin{bmatrix}
	1 \\ 0 \\ 0 \\ 0
\end{bmatrix}^\top,~
K = \begin{bmatrix}
1 \\ 1 \\ 1 \\ 1
\end{bmatrix}^\top $.
This system satisfies the condition of partial causal detectability with respect to $K$. Hence, it follows from Theorem \ref{thm:ODE:estimator} that there exists a functional ODE estimator of the form \eqref{obsv}.
	
We now design a functional ODE estimator for the given system by following the procedure provided in the proof of Theorem \ref{thm:ODE:estimator}. \\
\textbf{Step $1$:} By Lemma \ref{lm:Cdecomposition} and the (adaptation of the) algorithm provided in \cite{jaiswal2024existence} we obtain $U_O = I_4$,  $V_O = I_5$ and the following coefficient matrices for the reduced system:
\begin{eqnarray*}
 E_O = E,~A_O = A,~ B_O = B,~ C_O = C, \text{ and } K_O = K. 
\end{eqnarray*}
\textbf{Step $2$:} Using the method provided in \cite{berger2013addition}, we obtain the following matrices to convert the reduced system in QKF \eqref{eq:qkf}:
$P = \begin{bmatrix}
	0 & I_3 & 0 \\
	1 & 0 & 0 \\
	0 & 0 & 1
\end{bmatrix}$ and $Q = \begin{bmatrix}
	0 & 1 \\
	I_3& 0
\end{bmatrix}.$

\noindent \textbf{Step $3$:} This system does not contain positive finite eigenvalue and, hence
$U_1 = I_2,~J_f = \begin{bmatrix}
-1 & 1 \\ 0 & -1 \end{bmatrix}$.

\noindent \textbf{Step $4$:} $E_{\eta} = \begin{bmatrix}
1 & 0 \end{bmatrix}^\top$ is already in the required form, thus $U_2 = 1$. \\
\textbf{Step $5$:} Therefore, in the new coordinates the system becomes
\begin{eqnarray*}
\dot{x}_f(t) &=& \begin{bmatrix}
				-1 & 1 \\ 0 & -1
\end{bmatrix} x_f(t) + \begin{bmatrix}
				1 & 0 \\ 1 & 0
\end{bmatrix} \bar{u}(t) \\
0 &=& x_{\sigma}(t) + \begin{bmatrix}
1 & 0 \end{bmatrix} \bar{u}(t), \\
\dot{x}_{\eta}(t) &=& x_{\eta}(t) + \begin{bmatrix}
1 & 0 \end{bmatrix} \bar{u}(t), \\
0 &=& x_{\eta}(t) + \begin{bmatrix}
0 & -1 \end{bmatrix} \bar{u}(t), \\
z(t) &=& \begin{bmatrix}
      1 & 1
\end{bmatrix}x_f(t) + x_{\sigma}(t) + x_{\eta}(t) .
\end{eqnarray*}
\textbf{Step $6$:} This system has no $\epsilon$- and $f_1$-blocks. \\
\textbf{Step $7$:} From Step $5$, 
we obtain
\begin{eqnarray*}
x_{\sigma}(t) = \begin{bmatrix}
-1 & 0 \end{bmatrix} \bar{u}(t) \text{ and } x_{\eta}(t) = \begin{bmatrix}
0 & 1 
\end{bmatrix} \bar{u}(t). 
\end{eqnarray*}
\textbf{Step $8$:} Thus, in the new coordinates, the problem of functional ODE estimator design for the given system reduces to the problem of functional ODE estimator design for
\begin{eqnarray*}
\dot{x}_f(t) &=& \begin{bmatrix}
				-1 & 1 \\ 0 & -1
\end{bmatrix} x_f(t) + \begin{bmatrix}
				1 & 0 \\ 1 & 0
\end{bmatrix} \bar{u}(t), \\
z(t) &=& \begin{bmatrix}
				1 & 1
\end{bmatrix}x_f(t) + \begin{bmatrix}
				-1 & 1
\end{bmatrix}\bar{u}(t) .
\end{eqnarray*}
\textbf{Step $9$:} Since $x_{\eta}$ is obtained in Step $7$ above this step can be skipped. \\
\textbf{Step $10$:} Finally, we obtain the functional ODE estimator for the given system as follows:
\begin{eqnarray*}
\dot{w}(t) &=& \begin{bmatrix}
				-1 & 1 \\ 0 & -1
\end{bmatrix}w(t) + \begin{bmatrix}
				1 & 0 \\ 1 & 0
\end{bmatrix}\begin{bmatrix}
				u(t) \\ y(t)
\end{bmatrix} \\
\hat{z}(t) &=& \begin{bmatrix}
				1 & 1
\end{bmatrix}w(t) + \begin{bmatrix}
				-1 & 1
\end{bmatrix}\begin{bmatrix}
				u(t) \\ y(t)
\end{bmatrix}.
\end{eqnarray*}
Simulation results conducted in MATLAB are shown in Figures \ref{fig:fig1} and \ref{fig:fig2}. It can be observed that the proposed new design method provides an  asymptotic estimate $\hat z$ for the given functional $z$. In addition, it is clear from Figure \ref{fig:fig2} that the proposed functional ODE estimator is not a functional ODE observer, \emph{i.e.}, it does not exhibit the state matching property.

\begin{figure}[H]
	\centering
	\subfigure[Time response of state $z(t)$]{\includegraphics[width=.49\linewidth]{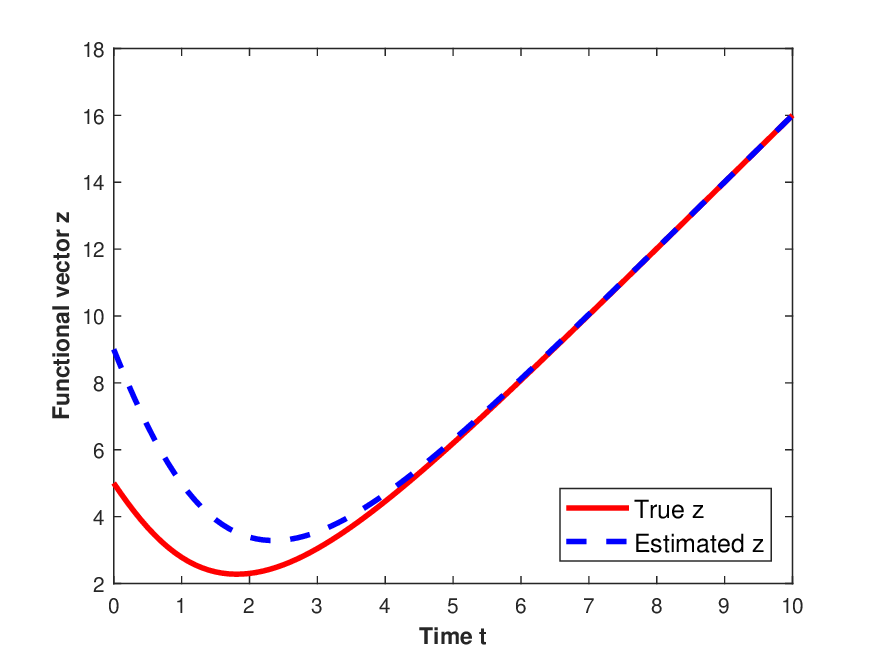}}
	\subfigure[Time response of estimation error in $z$]{\includegraphics[width=.49\linewidth]{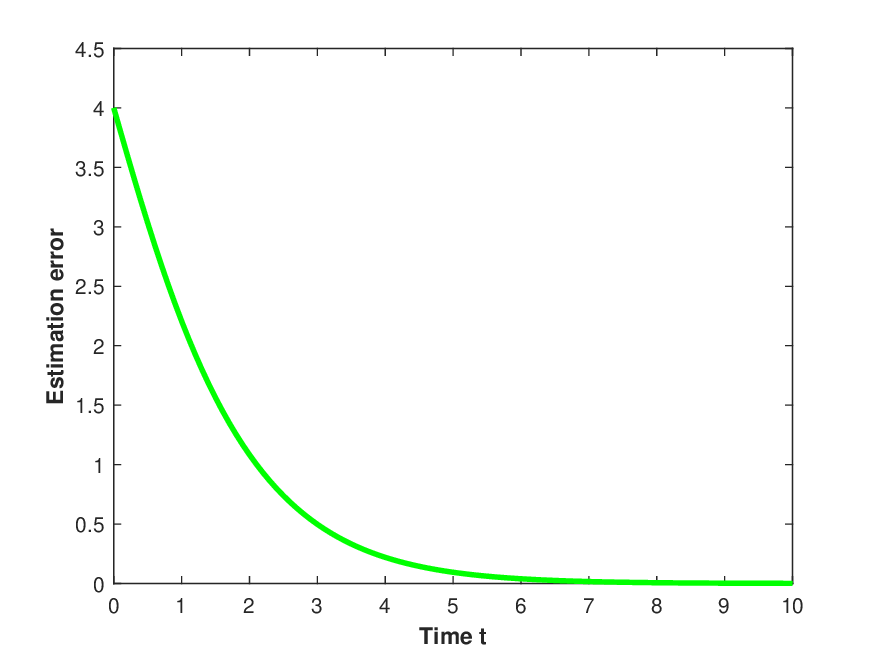}}
	\caption{Plot of estimated functional and estimation error with initial conditions $x(0) = \begin{bmatrix}
			1 & 2 & 3 & 0 \end{bmatrix}^T$, $w(0) =  \begin{bmatrix} 4 & 5 \end{bmatrix}$, and input $u(t) = t$.}
	\label{fig:fig1}
\end{figure}

\begin{figure}
	\centering
	\subfigure[Time response of state $z(t)$]{\includegraphics[width=.49\linewidth]{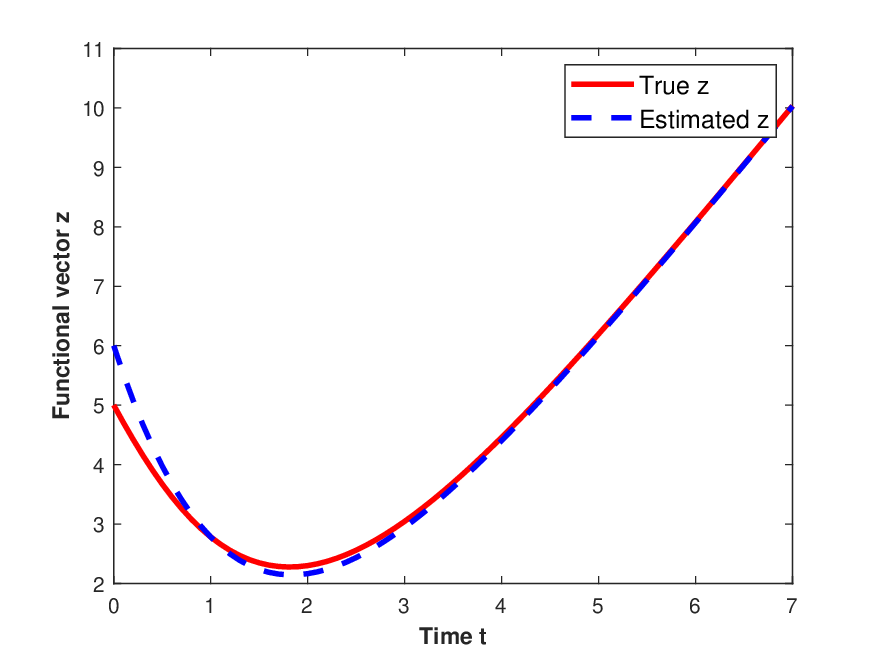}}
	\subfigure[Time response of estimation error in $z$]{\includegraphics[width=.49\linewidth]{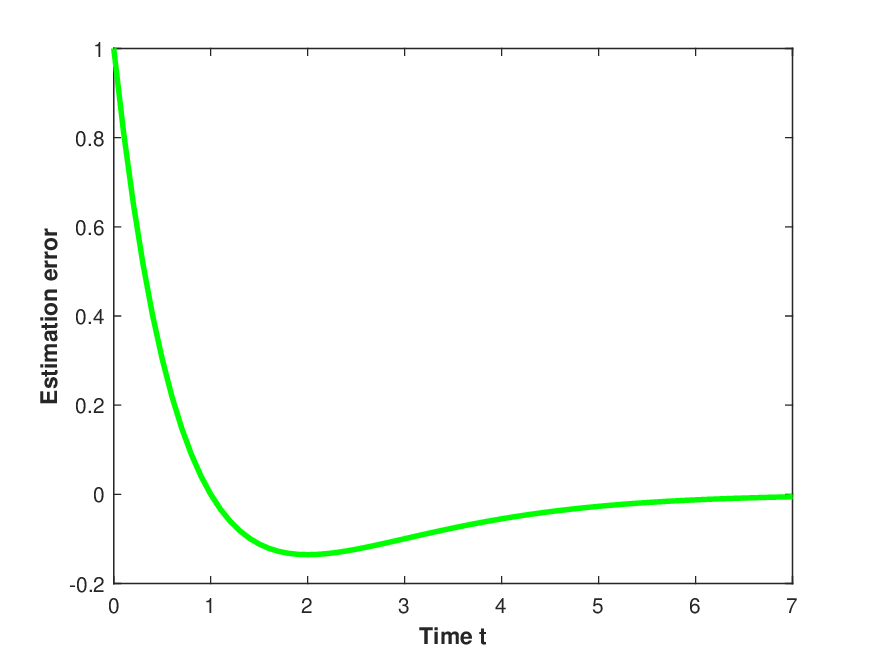}}
	\caption{Plot of estimated functional and estimation error with initial conditions $x(0) = \begin{bmatrix}
			1 & 2 & 3 & 0 \end{bmatrix}^T$, $w(0) =  \begin{bmatrix} 4 & 2 \end{bmatrix}$, and input $u(t) = t$}
	\label{fig:fig2}
\end{figure}
	
Now, we claim that there exists no functional ODE observer for this system, which is suggested by the fact that this system does not satisfy the existence conditions proposed in \cite{jaiswal2024existence,jaiswal2021functional,jaiswal2021existence,jaiswal2022functional}.
To see this, assume that  there exists a functional ODE observer of the form \eqref{obsv} and let $(x, 0, 0, Kx) \in \mathscr{B}$ be arbitrary for the given system. Then $w = 0$ and $\hat{z} = 0$ satisfy \eqref{obsv} with $u = 0$ and $y = 0$.  Since \eqref{obsv} is a functional ODE observer for the given system, we find that
	\begin{eqnarray*}
e(t) :=  z(t) - \hat{z}(t) = Kx(t)\to 0 \text{ for } t \to \infty
~\text{and } ~ e(0) = 0 \implies e(t) = 0, ~~\forall ~t >0.
\end{eqnarray*}
For instance, let us take the initial condition as $x(0) = \begin{bmatrix}
		0 & -1 & 1 & 0
\end{bmatrix}^\top$, then the solution of the system is $x(t) = \begin{bmatrix}
		0 & (t-1)e^t & e^t & 0
\end{bmatrix}^\top$, $y(t) = 0$, and $z(t) = te^t$ for $t \geq 0$. Here $e(0) = 0$ but $z(t) \neq 0 = \hat{z}(t)$ for all $t>0$. Thus, there exists no functional ODE observer for this system.
\end{example}

\section{Conclusion}\label{sec:conc}
A physically meaningful concept of partial causal detectability for LTI descriptor systems \eqref{dls} has been introduced, which is a natural extension of causal detectability of \eqref{dls} for $K = I_n$. Also, various equivalent characterizations of partial causal detectability have been established. Moreover, it has been proved that the notion of partial causal detectability is necessary and sufficient for the existence of functional ODE estimators. Remarks \ref{rem:obsv} and Example \ref{exp2} clarify that the concept of ODE observer and ODE estimator are not the same when $K \neq I_n$. Till date, the proposed existence condition in \cite{jaiswal2024existence} is the mildest known sufficient condition for the existence of a functional ODE observer. However, conditions which are necessary and sufficient for the existence of a functional ODE observer are not known. Future research directions include the development of some physical characterization to fill the gap between functional ODE observers and functional ODE estimators.

\bibliographystyle{siamplain}
\bibliography{Partial_causality}

\end{document}